\documentclass{article}
\usepackage[utf8]{inputenc}
\usepackage[margin=1in, includehead, includefoot]{geometry}
\usepackage{amsthm,amssymb,amsmath}
\usepackage{enumitem}
\usepackage{mathrsfs}
\usepackage{fancyhdr}
\usepackage{indentfirst}
\usepackage{graphicx}
\usepackage{placeins}
\usepackage{wrapfig}
\usepackage{mathtools}
\usepackage{tikz}
\usepackage{array,multirow}
\usepackage{algorithm,algpseudocode}
\usepackage{url}

\usepackage{xcolor}

\newtheorem{definition}{Definition}
\newtheorem{lemma}{Lemma}
\newtheorem{theorem}{Theorem}
\newtheorem{prop}{Proposition}
\newtheorem{corollary}{Corollary}
\newtheorem{observation}{Observation}

\DeclarePairedDelimiter\ceil{\lceil}{\rceil}

\newcommand{\toroman}[1]{\textit{\expandafter{\romannumeral #1\relax}}}

\newcommand{\cbeginproof}[0]{\par\noindent\textit{Proof.} }
\newcommand{\cendproof}[0]{ \qed\par\vspace{1em}}

\newcommand{\npcompleteproblem}[3]{ \par\vspace{0.5em}\noindent{\textbf{#1}}\newline \textbf{INSTANCE: } #2 \newline \textbf{QUESTION: } #3 \par\vspace{0.5em} }

\newcounter{algnum}

\title{On Redundant Locating-Dominating Sets}
\author{
    \small Devin C. Jean\\
    \small Electrical Engineering and \\ \small Computer Science Department \\
    \small Vanderbilt University\\
    \small \texttt{devin.c.jean@vanderbilt.edu}
    \and
    \vspace{0.6em}
    \small Suk J. Seo\\
    \vspace{0.6em}
    \small Computer Science Department\\
    \small Middle Tennessee State University\\
    \small \texttt{Suk.Seo@mtsu.edu}
}
\date{}

\begin{document}
\maketitle
\thispagestyle{empty}
\begin{abstract}
A \emph{locating-dominating} set in a graph $G$ is a subset of vertices representing ``detectors" which can locate an ``intruder" given that each detector covers its closed neighborhood and can distinguish its own location from its neighbors.
We explore a fault-tolerant variant of locating-dominating sets called \emph{redundant locating-dominating} sets, which can tolerate one detector malfunctioning (going offline or being removed).
In particular, we characterize redundant locating-dominating sets and prove that the problem of determining the minimum cardinality of a redundant locating-dominating set is NP-complete.
We also determine tight bounds for the minimum density of redundant locating-dominating sets in several classes of graphs including paths, cycles, ladders, $k$-ary trees, and the infinite hexagonal and triangular grids.
We find tight lower and upper bounds on the size of minimum redundant locating-dominating sets for all trees of order $n$, and characterize the family of trees which achieve these two extremal values, along with polynomial time algorithms to classify a tree as minimum extremal or not.
\end{abstract}

\noindent
\textbf{Keywords:} \textit{locating-dominating sets,  fault-tolerant, redundant locating-dominating sets, characterization, NP-complete, extremal trees, density}
\vspace{1em}

\noindent
\textbf{Mathematics Subject Classification:} 05C69

\section{Introduction}\label{sec:intro}
Let $G$ be a graph with vertices $V(G)$ and edges $E(G)$.

\begin{definition}
Given $\Re \subseteq \mathscr{P}(V(G))$, vertices $v,u \in V(G)$ are \emph{separable} if $\exists A \in \Re$ such that $v \in A \oplus u \in A$.
\end{definition}

\begin{definition}[\cite{old}]\label{def:disty-set}
$\Re \subseteq \mathscr{P}(V(G))$ is a \emph{distinguishing set} if every distinct pair of vertices is separable.
\end{definition}

In practice, a distinguishing set is often \emph{detector-based}, meaning that it is generated from a set of vertices representing the positions of detectors or sensors in the graph.
The method to generate a distinguishing set from a detector set varies depending on the capabilities of the detectors being used.
Let $S \subseteq V(G)$ be a set of detectors and $v \in S$.
We think of detector vertex $v$ as having one or more physical sensors, each of which may have a different \emph{detection region}: the area in which the physical sensor can sense an intruder.
Thus, $v$ is associated with a set of detection regions, denoted by $R(v) \subseteq \mathscr{P}(V(G))$.
Then the generated set $\cup_{v \in S}{R(v)}$ is a distinguishing set for sufficient choices of $S$.

The concept of a distinguishing set has also been studied under the name ``identifying system", introduced by Auger et al. \cite{watching-sys}.
Our introduction of distinguishing sets being generated from sets of detection regions $R(v) \subseteq \mathscr{P}(V(G))$ can be considered a generalization of their identifying systems.
In one approach, Karpovsky et al. \cite{karpovsky} gives each detector vertex $v$ a ``ball" region covering all vertices which are at most distance $r$ from $v$---which we will denote as $B_r(v) = \{w \in V(G) : d(v,w) \le r\}$---in which $v$ can sense intruders; this use of ball regions is a special case of our detection regions where $R(v) = \{B_r(v)\}$.
In another approach taken by Auger et al. \cite{watching-sys} each detector is assigned a ``watching zone", which is a subset of vertices that the detector can cover; this is closest to our set of detection regions $R(v)$, being a special case where $|R(v)| = 1$.
Many other approaches have also been taken; Lobstein \cite{dombib} maintains a bibliography of distinguishing-set-related parameters, which currently has over 440 papers.

\begin{definition}
The \emph{open-neighborhood} of a vertex $v \in V(G)$, denoted $N(v)$, is the set of all vertices adjacent to $v$: $\{w \in V(G) : vw \in E(G)\}$.
\end{definition}

\begin{definition}
The \emph{closed-neighborhood} of a vertex $v \in V(G)$, denoted $N[v]$, is the set of all vertices adjacent to $v$, as well as $v$ itself: $N(v) \cup \{v\}$.
\end{definition}

The critical difference between different types of \emph{detector-based distinguishing (DBD)} sets is the choice of detection regions, $R$.
For instance, \emph{identifying codes (ICs)} are DBD sets with $R(v) = \{N[v]\}$, \emph{open-locating-dominating (OLD)} sets are DBD sets with $R(v) = \{N(v)\}$, and \emph{locating-dominating (LD)} sets are DBD sets with $R(v) = \{\{v\}, N(v)\}$.

\begin{definition}\label{def:general-dom}
For a DBD set $S \subseteq V(G)$, a vertex $v \in V(G)$ is \emph{$k$-dominated} if $|\{ u \in S : v \in \cup{R(u)} \}| = k$.
\end{definition}
\begin{definition}
The \emph{domination count} of a vertex $v$, denoted $dom(v)$, is $k$ if and only if $v$ is $k$-dominated.
\end{definition}

An important real-world application of distinguishing sets is in the creation of automated security systems for locating an intruder in a facility, for locating a faulty processor in a multiprocessor network, etc.; in any case, we term the phenomenon being detected the ``intruder".
The information required to locate the intruder would come from physical detectors, hence the utility of DBD sets.
Let $S \subseteq V(G)$ be the DBD set with detection regions $R$.
In typical applications, a detector vertex $v \in S$ transmits some unique signal for every distinct intersection of elements in $R(v)$; that is, for every possible combination of overlapping detection regions, of which there are precisely $c$ distinct intersections, where $c = |\{\cap{A} : A \subseteq R(v)\}|$.
Note that this calculation of $c$ includes the empty set, which denotes that no intruder was detected by the sensor.
In many applications, the elements of $R(v)$ are disjoint, in which case $c = |R(v)| + 1$.
Typically, we think of the signals as transmitting an integral value in $[0, c)$, where 0 denotes no intruder being detected.
From the raw transmitted information, the system applies separability to either determine the exact location of the intruder or conclude that there is none.
We will assume that, at any given time, there is at most one intruder.

In this paper, we primarily focus on LD sets, which use $R(v) = \{\{v\},N(v)\}$.
For convenience, we define a transmitted value of 0 to indicate no intruder being detected, 1 to indicate an intruder in the open neighborhood $N(v)$, and 2 to indicate an intruder at its location $\{v\}$.
This assignment of 0, 1, and 2 is the same convention as used in describing LD sets when first introduced by Slater \cite{ftld,dom-ref-sets,dom-loc-acyclic}.

\begin{definition}[\cite{dom-ref-sets}]
An LD set $S \subseteq V(G)$ is a subset of vertices such that $\forall v,u \in V(G)-S$ with $v \neq u$, $\varnothing \neq N(v) \cap S \neq N(u) \cap S$.
\end{definition}

As the application of distinguishing sets is in the real world, it is often the case that a detector in the network may be faulty.
These type of errors can be modeled by fault-tolerant variants of LD sets, which include detector redundancies and the system's ability to handle false negatives and false positives.
In this paper, we will explore \emph{redundant location-dominating (RED:LD)} sets, which can tolerate one detector malfunctioning or being removed.
More general types of detector-based fault-tolerance have also been studied by Seo and Slater \cite{ftsets, gen}.

\begin{definition}\label{def:red-ld}
A \emph{redundant LD (RED:LD)} set is an LD set $S \subseteq V(G)$ such that for any detector $ v \in S$, $S-\{v\}$ is also an LD set.
\end{definition}

\begin{wrapfigure}{r}{0.24\textwidth}
    \centering
    \begin{tabular}{c}
        \includegraphics[width=0.17\textwidth]{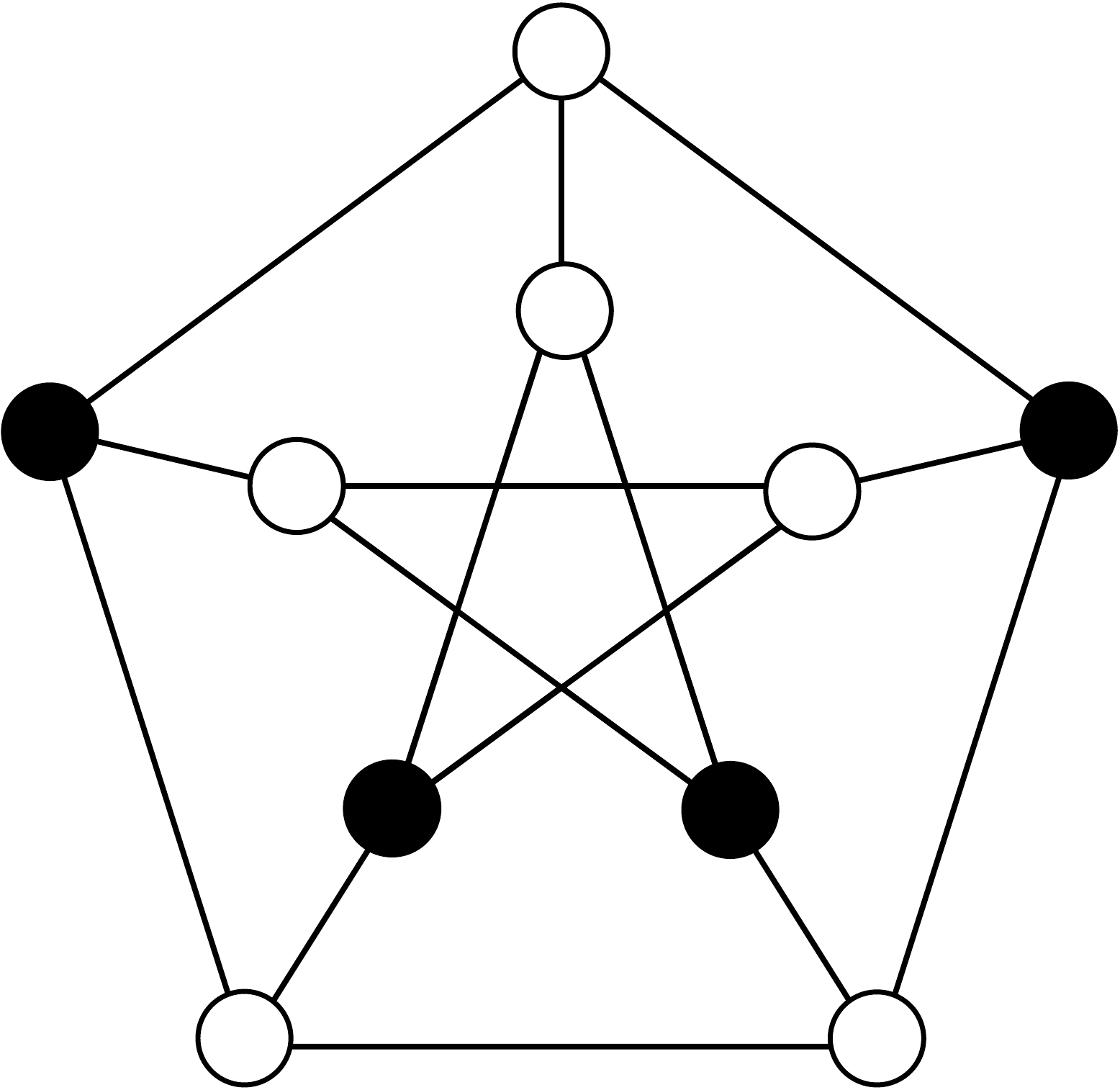} \\ (a) \\\\ \includegraphics[width=0.17\textwidth]{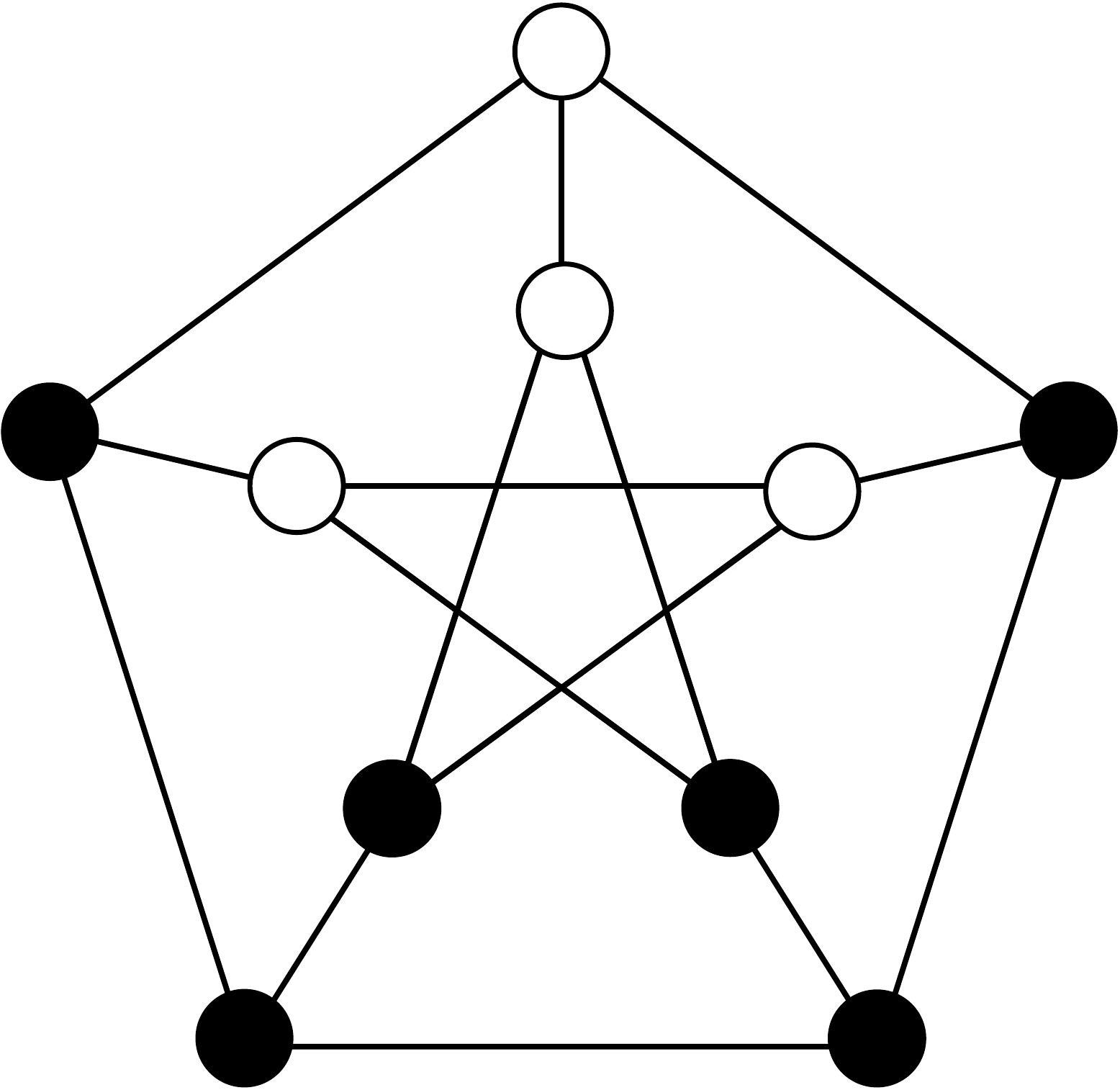} \\
        (b)
    \end{tabular}
    \caption{Optimal LD (a) and RED:LD (b) sets on the Petersen graph.}
    \label{fig:petersen-examples}
\end{wrapfigure}

Any superset of a DBD set is clearly also a DBD set, so we are interested in the smallest sets with the given properties; this is especially important in real-world applications, as each detector represents a piece of physical hardware, making the smallest DBD set the most cost-effective.
For a finite graph $G$, let $\textrm{LD}(G)$ and $\textrm{RED:LD}(G)$ denote the cardinality of the smallest such sets in $G$.
For infinite graphs, we measure via the \textit{density}, which is defined to be the ratio of the number of detectors to $|V(G)|$.
The minimum density of such a set in $G$ is denoted by $\textrm{LD}\%(G)$ and $\textrm{RED:LD}\%(G)$.
In some cases, we may prefer to use densities for finite graphs instead of cardinality.


As an example, consider the constructions of LD and RED:LD sets on the Petersen graph, $G$, as shown in Figure~\ref{fig:petersen-examples}.
One can verify that the sets of shaded vertices satisfy the requirements of their corresponding definitions, and so $\textrm{LD}(G) \le 4$ and $\textrm{RED:LD}(G) \le 6$.
Additionally, no smaller set exists which satisfies the requirements; therefore, these are optimal, and we have $\textrm{LD}(G) = 4$ and $\textrm{RED:LD}(G) = 6$.
If we would prefer to use densities, we have $\textrm{LD}\%(G) = \frac{4}{10} = \frac{2}{5}$ and $\textrm{RED:LD}\%(G) = \frac{6}{10} = \frac{3}{5}$.

In the following section, we characterize RED:LD sets on arbitrary graphs and establish the existence criteria.
In Section~\ref{sec:np-complete} we prove the NP-completeness of the problem of determining the minimum cardinality of a RED:LD set.
In Section~\ref{sec:graphs} we explore finding RED:LD sets in several special classes of graph, including cycles, ladders, trees, and some infinite grids.
Additionally, we find tight lower and upper bounds on the size of minimum RED:LD sets for all trees of order $n$, and characterize the trees which achieve these two extremal values.

\FloatBarrier
\section{Characterization and Existence}\label{sec:char}

We will begin by establishing necessary and sufficient requirements for a RED:LD set in an arbitrary graph, $G$.
A characterization of (plain) LD sets is given by Theorem~\ref{theo:ld-char} for comparison.

\begin{theorem}[\cite{dom-ref-sets}]\label{theo:ld-char}
A detector set $S \subseteq V(G)$ is an LD set if and only if the following are true:
\begin{enumerate}[noitemsep, label=\roman*.]
    \item $\forall v \in V(G)-S$, $|N(v) \cap S| \ge 1$
    \item $\forall v,u \in V(G)-S$ with $v \neq u$, $|(N(v) \cap S) \triangle (N(u) \cap S)| \ge 1$
\end{enumerate}
\end{theorem}

Blidia et al. \cite{ld-tree-char} proved a specific characterization for unique minimum LD sets in trees.
Hernando et al. \cite{locdom-bounds} found Nordhaus-Gaddum bounds---tight bounds for a graphical parameter on sums and products of $G$ with $\overline{G}$---for LD sets on families of graphs, and characterized the families.
To the best of our knowledge, there exists no previous characterization of RED:LD.

\begin{lemma}\label{lem:red-ld-s-imp}
Let $S \subseteq V(G)$ be a RED:LD set and $v \in S$. Then
\begin{enumerate}[noitemsep, label=\roman*.]
    \item $|N(v) \cap S| \ge 1$
    \item $\forall u \in V(G)-S$, $|((N(v) \cap S) \triangle (N(u) \cap S)) - \{v\}| \ge 1$
\end{enumerate}
\end{lemma}
\begin{proof}
Suppose property \toroman{1} is false; then $\exists v \in S$ such that $|(N(v) \cap S)| = 0$.
Because $v$ is the only detector that can detect an intruder at $v$, the set $S - \{v\}$ cannot find an intruder at $v$, a contradiction.
Next, suppose that property \toroman{2} is false.
Then $\exists v \in S$ and $\exists u \in V(G)-S$ such that $((N(v) \cap S) \triangle (N(u) \cap S)) - \{v\} = \varnothing$.
Note that $v \notin N(v)$, so $v \in (N(v) \cap S) \triangle (N(u) \cap S)$ if and only if $v \in N(u)$.
If $v \in N(u)$, then $(N(v) \cap S) \triangle (N(u) \cap S) = \{v\}$, otherwise $v \notin N(u)$ and $(N(v) \cap S) \triangle (N(u) \cap S) = \varnothing$.
In either case, we see the new set $S'=S-\{v\}$ results in $v,u \notin S'$ but $(N(v) \cap S') \triangle (N(u) \cap S') = \varnothing$, contradicting that $S'$ is an LD set, completing the proof.
\end{proof}

\begin{lemma}\label{lem:red-ld-not-s-imp}
Let $S \subseteq V(G)$ be a RED:LD set and $v \notin S$. Then
\begin{enumerate}[noitemsep, label=\roman*.]
    \item $|N(v) \cap S| \ge 2$
    \item $\forall u \in V(G)-S$ with $u \neq v$, $|(N(v) \cap S) \triangle (N(u) \cap S)| \ge 2$
\end{enumerate}
\end{lemma}
\begin{proof}
Suppose that property \toroman{1} is false; then $\exists v \in V(G)-S$ such that $|N(v) \cap S| \le 1$.
Additionally, note that because $v \notin S$, we require $|N(v) \cap S| \ge 1$ in order for $S$ to be an LD set; so we assume $|N(v) \cap S| = 1$.
Then $N(v) \cap S = \{t\}$ for some $t \in S$.
The set $S'=S-\{t\}$ cannot detect an intruder at $v$, a contradiction.
Next, suppose that property \toroman{2} is false; then $\exists v,u \in V(G)-S$ with $v \neq u$ such that $|(N(v) \cap S) \triangle (N(u) \cap S)| \le 1$.
Note that $|(N(v) \cap S) \triangle (N(u) \cap S)| \neq 0$ because this would contradict that $S$ is an LD set; therefore, we assume that $|(N(v) \cap S) \triangle (N(u) \cap S)| = 1$, meaning $(N(v) \cap S) \triangle (N(u) \cap S) = \{t\}$ for some $t \in S$.
Consider the set $S'=S-\{t\}$; then $v,u \notin S'$ but $N(v) \cap S' = N(u) \cap S'$, a contradiction.
\end{proof}

\begin{theorem}\label{theo:red-ld-char}
A set $S \subseteq V(G)$ is a RED:LD set if and only if the following are true:
\begin{enumerate}[noitemsep, label=\roman*.]
    \item $\forall v \in V(G)$, $|N[v] \cap S| \ge 2$
    \item $\forall v \in S$ and $\forall u \in V(G)-S$, $|((N(v) \cap S) \triangle (N(u) \cap S)) - \{v\}| \ge 1$
    \item $\forall v,u \in V(G)-S$ with $u \neq v$, $|(N(v) \cap S) \triangle (N(u) \cap S)| \ge 2$
\end{enumerate}
\end{theorem}
\begin{proof}
If $S$ is a RED:LD set, then properties \toroman{1}--\toroman{3} are given by Lemmas \ref{lem:red-ld-s-imp} and \ref{lem:red-ld-not-s-imp}.
For the converse, suppose $S \subseteq V(G)$ satisfies properties \toroman{1}--\toroman{3}.
Properties \toroman{1} and \toroman{3} together are sufficient to invoke Theorem~\ref{theo:ld-char}; thus, $S$ is an LD set.
Now, suppose we remove a detector $v \in S$, creating a new set $S'=S-\{v\}$.
Let $a \in V(G)-S' = (V(G)-S) \cup \{v\}$.
If $a \in V(G)-S$, then property \toroman{1} gives us that $|N(a) \cap S| \ge 2$, which means $|N(a) \cap S'| \ge 1$.
Otherwise, $a = v \in S$; therefore, property \toroman{1} yields that $|N(a) \cap S| \ge 1$, and because $v \notin N(v)$ we know $|N(a) \cap S'| \ge 1$.
Thus, $a$ is at least 1-open-dominated by $S'$.
Next, let $b \in V(G)-S'$ with $a \neq b$.
If $a \in V(G)-S$ and $b \in V(G)-S$, then property \toroman{3} yields that $|(N(v) \cap S) \triangle (N(u) \cap S)| \ge 2$.
This means that the open-neighborhoods intersected with $S$ have at least two differences; in removing a single detector $v \in S$, we eliminate at most one of the differences, so $|(N(v) \cap S') \triangle (N(u) \cap S')| \ge 1$.
Otherwise, without loss of generality, let $a = v \in S$, which requires $b \in V(G)-S$ because $b \in (V(G)-S) \cup \{v\}$ and $a \neq b$ by hypothesis.
If $a \in N(b)$, then property \toroman{2} yields that $|(N(a) \cap S) \triangle (N(b) \cap S)| \ge 2$, meaning $|(N(a) \cap S') \triangle (N(b) \cap S')| \ge 1$.
Otherwise, $a \notin N(b)$ and property \toroman{2} gives us that $|(N(a) \cap S) \triangle (N(b) \cap S)| \ge 1$; in this case $N(a) \cap S = N(a) \cap S'$ and $N(b) \cap S = N(b) \cap S'$, so $|(N(a) \cap S') \triangle (N(b) \cap S')| \ge 1$.
Thus, $a$ and $b$ are 1-distinguished.
As we've now demonstrated that all vertices in $V(G)-S'$ are at least 1-distinguished and 1-dominated from one another, Theorem~\ref{theo:ld-char} yields that $S'$ is an LD set; and because $v \in S$ was chosen arbitrarily, $S$ is a RED:LD set.
\end{proof}

\begin{definition}\label{def:red-ld-disty}
For a RED:LD set $S \subseteq V(G)$ and $u,v \in V(G)$, $v$ is \emph{$k$-distinguished} from $u$ if $|((N(v) \cap S) \triangle (N(u) \cap S)) - \{v, u\}| \ge k$.
\end{definition}

\begin{definition}\cite{twin-free} 
Two distinct vertices $v,u \in V(G)$ are said to be \emph{twins} if $N[u] = N[v]$ (\emph{closed twins}) or $N(u) = N(v)$ (\emph{open twins}).
\end{definition}

With Definitions \ref{def:general-dom} and \ref{def:red-ld-disty}, we see that Theorem~\ref{theo:red-ld-char} requires every vertex be at least 2-dominated, that each detector/non-detector pair be 1-distinguished, and that each non-detector pair be 2-distinguished.
Clearly, a RED:LD set $S$ exists if and only if $\delta(G) \neq 0$, as detector vertices have no requirement other than being 2-dominated.
We also see that if $u$ and $v$ are twins, then we require $\{u,v\} \subseteq S$ in order to be distinguished.

\begin{observation}\label{obs:red-ld-exist}
A RED:LD set exists if and only if $\delta(G) \neq 0$.
\end{observation}

\begin{observation}\label{obs:red-ld-upper}
For a complete $k$-partite graph, if $k = 2$ or no part is a singleton, then $\textrm{RED:LD}(G) = n$.
\end{observation}


By Observation~\ref{obs:red-ld-exist}, we are guaranteed that a RED:LD set exists on any connected graph of order $n \ge 2$.
From Observation~\ref{obs:red-ld-upper}, we see that all complete graphs and bipartite complete graphs (including stars) have $\textrm{RED:LD}(G) = n$.

\vspace{0.8em}
By Theorem~\ref{theo:ld-char}, if a graph, $G$, has $\textrm{LD}(G) = k$, then there can be at most $2^k-1$ non-detectors, one for each non-empty subset of the $k$ detectors.

\begin{observation}\label{theo:max-v-ld}
If $LD(G) \le k$, then $|V(G)| \le 2^k + k - 1$.
\end{observation}

\begin{theorem}\label{theo:max-v-redld}
If $\textrm{RED:LD}(G) \le k$, then $|V(G)| \le 2^{k-1}+k-2$.
\end{theorem}
\begin{proof}
Suppose we have a RED:LD set, $S \subseteq V(G)$, with $|S| \le k$; then by definition there exists a LD set $S'$ with $|S'| \le k-1$.
Observation~\ref{theo:max-v-ld} gives us that $|V(G)| \le 2^{k-1}+k-2$.
\end{proof}

\begin{theorem}
If $k=2j$, there is a graph of size $2^{k-1} + k - 2$ with $\textrm{RED:LD}(G) = k$.
\end{theorem}
\begin{proof}
We begin with a complete graph on $k$ vertices, where every vertex is a detector.
We then add an additional $\binom{k}{2}$ non-detectors which are adjacent to a distinct pair of detectors, an additional $\binom{k}{4}$ non-detectors which are adjacent to distinct sets of 4 detectors, an additional $\binom{k}{6}$ non-detectors which are adjacent to distinct sets of 6 detectors, and so on through $\binom{k}{k-2}$ (as $k$ detectors representing $\binom{k}{k}$ were already created at the beginning).
It is easy to verify that ever vertex is at least 2-dominated, as each subset of the $k$ detectors was at least size 2, and all vertices are at least 2-distinguished because only even sized subsets were chosen.
Thus, we have $\left(\binom{k}{2} + \binom{k}{4} + \hdots + \binom{k}{k-2}\right) + k$ vertices in total.
The summation $\binom{k}{0}+\binom{k}{2}+\binom{k}{4}+\cdots+\binom{k}{k}$ is known to be $2^{k-1}$; thus, $|V(G)| = 2^{k-1} - 2 + k$, completing the proof.
\end{proof}

For $k=2j+1$, the largest $|V(G)|$ with $\textrm{RED:LD}(G) = k$ we have constructed has $|V(G)|= 2^{k-1} + \frac{k-5}{2}$.

\begin{theorem}\label{theo:max-redld}
A (connected) graph, $G$, with $n \ge 3$ has $\textrm{RED:LD}(G) = n$ if and only if every vertex is a leaf vertex, support vertex, or is a twin with some other vertex.
\end{theorem}
\begin{proof}
Firstly, from Theorem~\ref{theo:red-ld-char}, we see that all vertices must be at least 2-dominated, implying all leaf and support vertices must be detectors.
We also see that if $u,v \in V(G)$ are twins, then they must both be detectors in order to be distinguished.
Thus, if every vertex is a leaf, support, or twin vertex, $\textrm{RED:LD}(G) = n$.
For the converse, suppose that $v \in V(G)$ is a non-leaf, non-support vertex which is not a twin with any other vertex; let $S = V(G)-\{v\}$, and let $u \in S$.
Because $v$ is not a leaf vertex and it is the only non-detector, it is at least 2-dominated,
Additionally, because $v$ is not a leaf or support vertex, if $u$ is a leaf node, then it is 2-dominated; otherwise, $deg(u) \ge 2$, so $u$ is at least 2-dominated by itself and one or more of its neighbors.
Therefore, all vertices are at least 2-dominated.
From Theorem~\ref{theo:red-ld-char}, we see that two detectors have no distinguishing requirements, and there are no distinct pairs of non-detectors, so the only remaining requirement is showing that $v$ and our arbitrary $u$ are distinguished.
By hypothesis, $v$ is not a twin with any other vertex, so $N(u) \neq N(v)$ and $N[u] \neq N[v]$.
If $uv \notin E(G)$, then $((N(v) \cap S) \triangle (N(u) \cap S)) - \{u,v\} = N(v) \triangle N(u) \neq \varnothing$, so $u$ and $v$ are distinguished.
Otherwise, $uv \in E(G)$, then $((N(v) \cap S) \triangle (N(u) \cap S)) - \{u,v\} = (N(v) \triangle (N(u) - \{v\})) - \{u,v\} = N[v] \triangle N[u] \neq \varnothing$, so $u$ and $v$ are distinguished.
Thus, $S = V(G) - \{v\}$ is a RED:LD set for $G$, implying $\textrm{RED:LD}(G) < n$, completing the proof.
\end{proof}

From Theorems~\ref{theo:max-v-redld} and \ref{theo:max-redld}, we have the following corollary.

\begin{corollary}
Let $G$ be a graph; then $\frac{k}{2^{k-1} + k - 2} \leq \textrm{RED:LD\%}(G) \leq 1$.
\end{corollary}

\section{NP-Completeness}\label{sec:np-complete}

The problem of finding the value of LD($G$) for an arbitrary graph has been known to be NP-complete \cite{ld-np-complete-2,NP-complete-ld}---see \cite{np-complete-bible} for more information on NP-completeness.
We will prove that the problem of finding the value for RED:LD($G$) is also NP-complete.

\npcompleteproblem{3-SAT}{Let $X$ be a set of $N$ variables.
Let $\psi$ be a conjunction of $M$ clauses, where each clause is a disjunction of three literals from distinct variables of $X$.}{Is there is an assignment of values to $X$ such that $\psi$ is true?}

\npcompleteproblem{Redundant Locating-Domination (RED-LD)}{A graph $G$ and integer $K$ with $2 \le K \le |V(G)|$.}{Is there exists a RED:LD set $S$ with $|S| \le K$? Or equivalently, is RED:LD($G$) $\le K$?}

\begin{figure}[ht]
    \centering
    \includegraphics[width=0.7\textwidth]{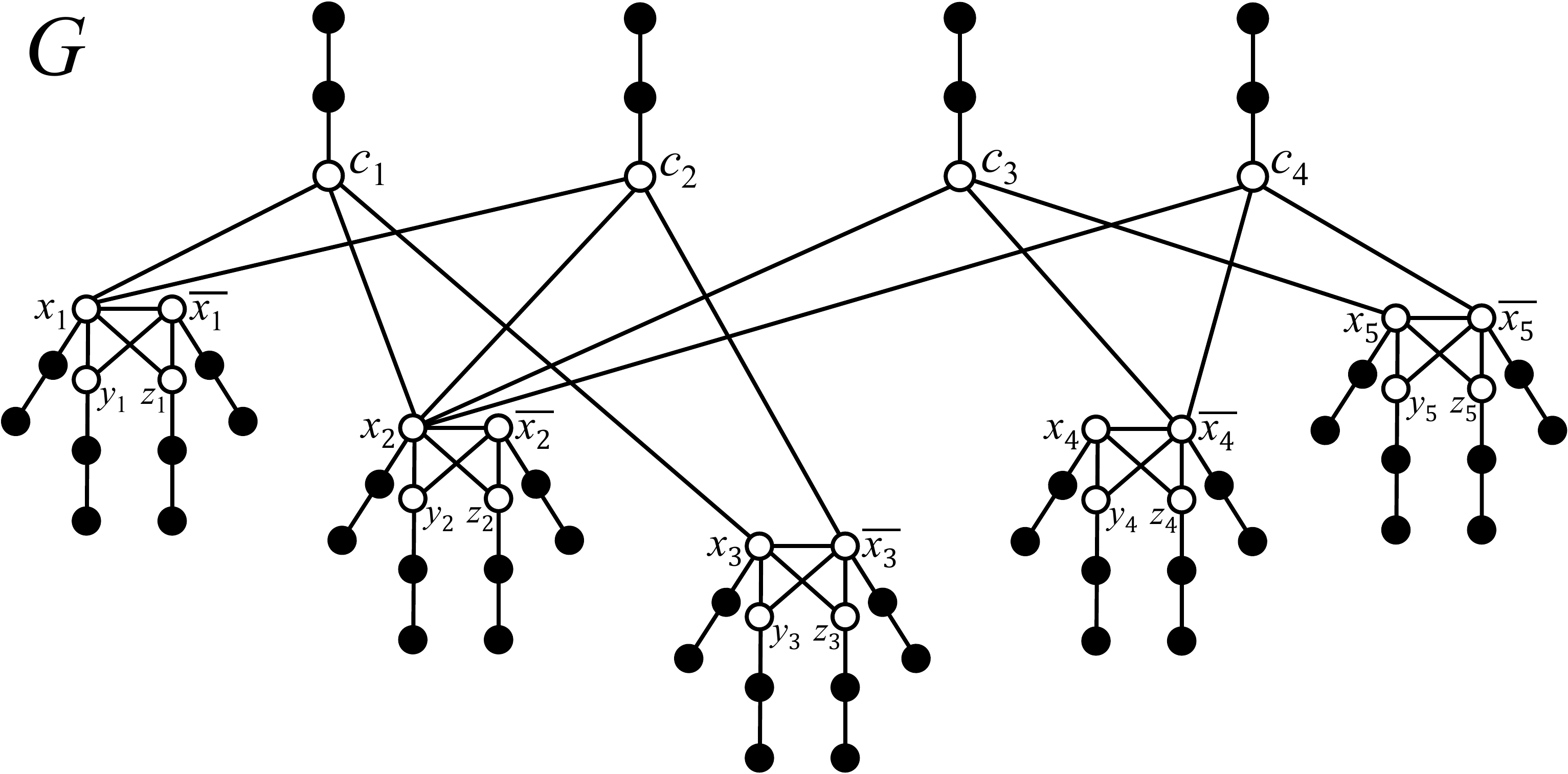}
    \caption{Construction of $G$ from $(x_1 \lor x_2 \lor x_3) \land (x_1 \lor x_2 \lor \overline{x_3}) \land (x_2 \lor \overline{x_4} \lor x_5) \land (x_2 \lor \overline{x_4} \lor \overline{x_5})$}
    \label{fig:example-clause}
\end{figure}

\begin{theorem}\label{theo:red-np-complete}
The RED-LD problem is NP-complete.
\end{theorem}
\cbeginproof
RED-LD is NP, as every possible candidate solution can be generated non-deterministically in polynomial time, and each candidate can be verified in polynomial time.
We will show a reduction from 3-SAT to RED-LD.

\begin{wrapfigure}[12]{r}{0.2\textwidth}
    \centering
    \begin{tabular}{c}
        \includegraphics[width=0.15\textwidth]{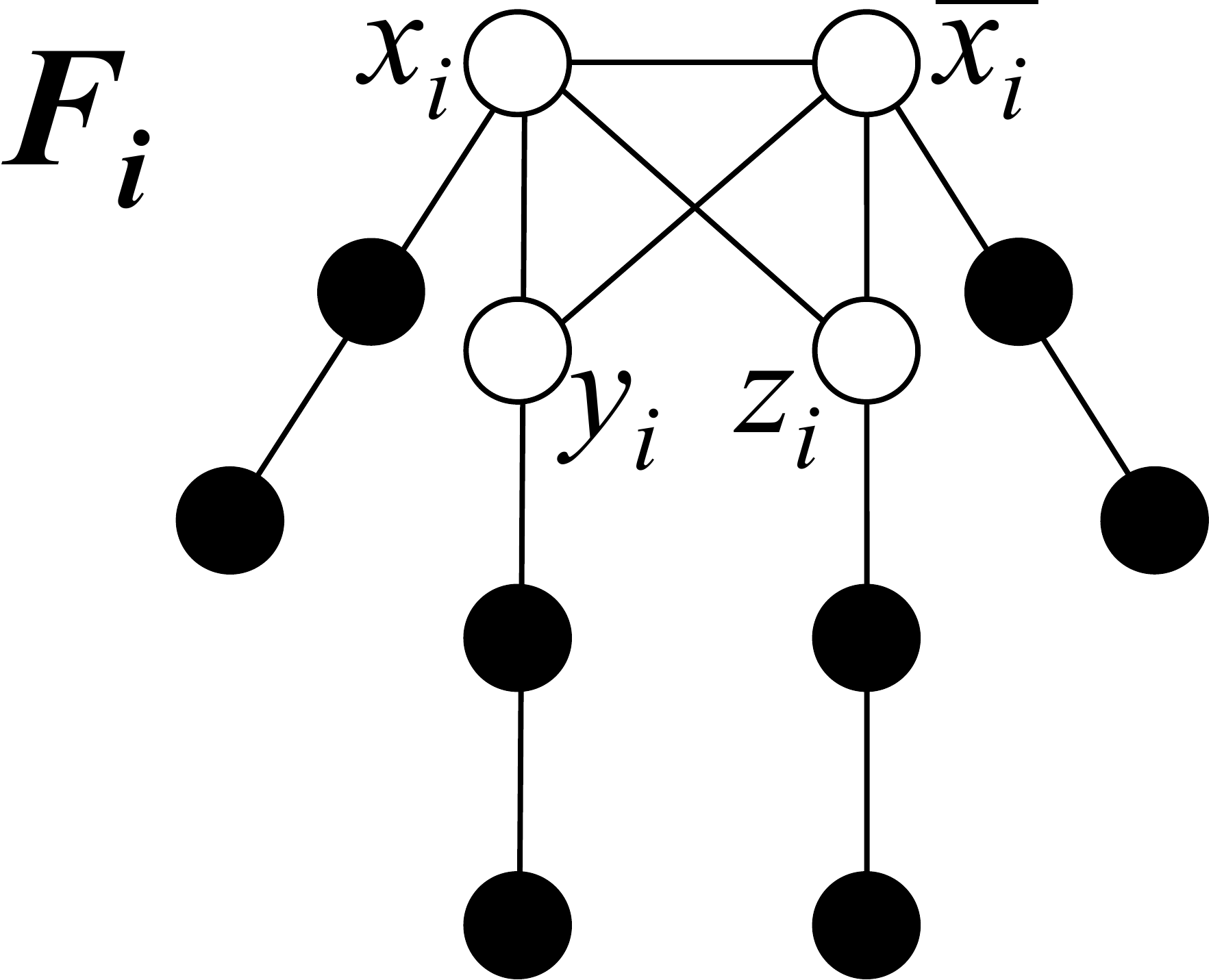} \\ \\ \\
        \includegraphics[width=0.075\textwidth]{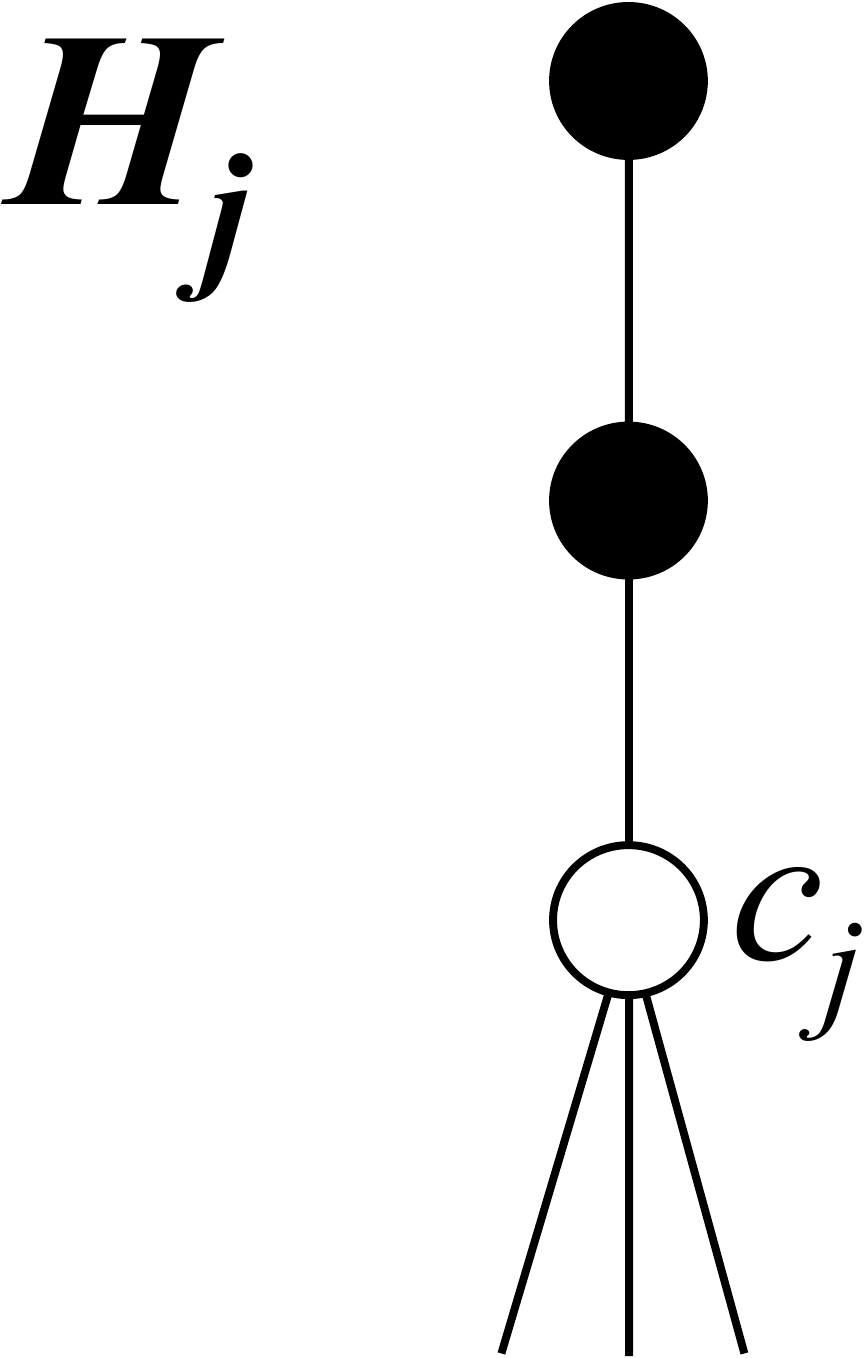}
    \end{tabular}
    \caption{Variable and Clause graphs.}
    \label{fig:variable-clause}
\end{wrapfigure}

Let $\psi$ be an arbitrary instance of the 3-SAT problem with $M$ clauses on $N$ variables.
We will construct a graph, $G$, as follows.
For each variable $x_i$, we create an instance of the $F_i$ graph, depicted in Figure~\ref{fig:variable-clause}~(a); note that there is a vertex for the variable, $x_i$, and its negation, $\overline{x_i}$.
For each clause $c_j$ of $\psi$, we create an instance of the $H_j$ graph, depicted in Figure~\ref{fig:variable-clause}~(b).
For each clause $c_j = \alpha \lor \beta \lor \gamma$, we create an edge from the $c_j$ vertex to $\alpha$, $\beta$, and $\gamma$ from the variable graphs, each of which is either some $x_i$ or $\overline{x_i}$; for an example, see Figure~\ref{fig:example-clause}.
The resulting graph has precisely $12N + 3M$ vertices and $13N + 5M$ edges, and can be constructed in polynomial time.

Suppose $S \subseteq V(G)$ is an optimal (minimum) RED:LD set on $G$.
By Theorem~\ref{theo:red-ld-char}, every vertex must be 2-dominated; thus, we require $8N + 2M$ detectors, as shown by the shaded vertices in Figure~\ref{fig:variable-clause}.
For this specific graph, $G$, it is the case that 2-domination of every vertex will be sufficient to distinguish every pair of vertices, as required by Theorem~\ref{theo:red-ld-char}.
For any $F_i$, we require $\{x_i,\overline{x_i},y_i\} \cap S \neq \varnothing$ to 2-dominate $y_i$ and $\{x_i,\overline{x_i},z_i\} \cap S \neq \varnothing$ to 2-dominate $z_i$; because $S$ is assumed to be a minimum set, it must be the case that $|\{x_i,\overline{x_i}\} \cap S| = 1$, as taking $y_i$ or $z_i$ would require two detectors.
Thus, $|S| \ge 9N + 2M$; if $|S| = 9N + 2M$, then $c_j \notin S$ for any $j$ and we see that $c_j$ must be 2-dominated by one of the three $x_i$ or $\overline{x_i}$ vertices it is adjacent to.
Therefore, each clause $c_j$ is true, and we have a satisfying truth assignment for the 3-SAT problem.

Now suppose we have a satisfying truth assignment to the 3-SAT problem.
For each variable, $x_i$, if $x_i$ is true then we let the vertex $x_i \in S$; otherwise, we let $\overline{x_i} \in S$.
By construction, this 2-dominates every vertex in $G$; thus $S$ is an optimal RED:LD set on $G$.
\cendproof

\FloatBarrier
\section{Special classes of graph}\label{sec:graphs}


\begin{wrapfigure}{r}{0.20\textwidth}
    \centering
    \includegraphics[width=0.20\textwidth]{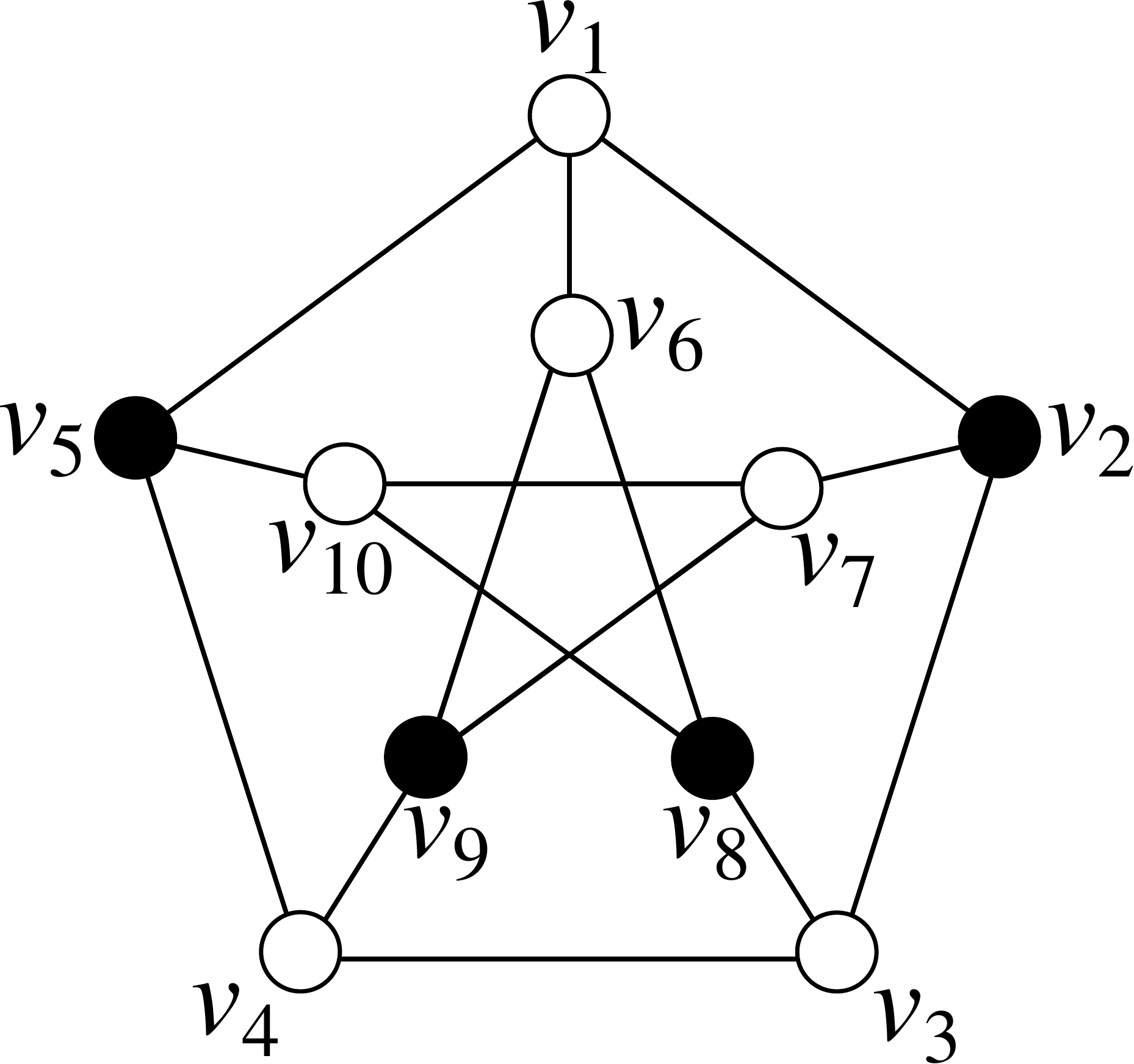}
    \caption{Example graph for share.}
    \label{fig:sh-ex}
\end{wrapfigure}

When establishing a lower bound for the density of any type of dominating set, it is convenient to use a \emph{share argument} \cite{ourtri,ourking,oldking}, first introduced by Slater \cite{ftld}.
In a share argument, instead of finding a lower bound for density directly, we invert the problem and find an upper bound for the maximum amount of sharing of the domination of each vertex.

For example, in  an LD set $S \subseteq V(G)$, a vertex $v \in V(G)$ is dominated by every detector in its closed neighborhood; this means $v$ is dominated $|N[v] \cap S|$ times.
Every dominator of $v$ has exactly $|N[v] \cap S|^{-1}$ of the share in dominating $v$; this is known as the \emph{partial share} of $v$, denoted $sh[v]$.
We group these partial shares in terms of the dominator, $x \in S$.
The \emph{share} of $x$, denoted $sh(x)$, is the sum of partial shares of every vertex dominated by $x$: $\sum_{w \in N[x]}{sh[w]}$.
By this construction, the average share of all detector vertices is equal to the inverse of the density of $S$ in $V(G)$.
Thus, we can invert the upper bound of average share to obtain a lower bound on the density.
For convenience, we will often shorten share sums using a ``sigma" notation, $\sigma_{A}$, which we define as $\sum_{a \in A}{\frac{1}{a}}$ where $A$ is a sequence of single-character symbols or numbers; thus, $\sigma_{2234} = \frac{1}{2} + \frac{1}{2} + \frac{1}{3} + \frac{1}{4}$.

Consider the Peterson graph $G$ from Figure~\ref{fig:sh-ex}, where detector set $S \subseteq V(G)$ is the set of shaded vertices.
Detector vertex $v_5 \in S$ dominates four vertices: $v_1$, $v_4$, $v_5$, and $v_{10}$.
Vertex $v_5$ is dominated only by $v_5$ (itself), while vertices $v_1$, $v_4$, and $v_{10}$ are dominated twice, by $v_5$ and some other detector.
Thus, $sh(v_5) = sh[v_1] + sh[v_4] + sh[v_5] + sh[v_{10}] = \sigma_{2212} =\frac{1}{2} + \frac{1}{2} + \frac{1}{1} + \frac{1}{2} = \frac{5}{2}$.
One can verify that detectors $v_2$, $v_8$, and $v_9$ also have a share of $\frac{5}{2}$ by applying similar logic.
Thus, the average share of all detectors is $\frac{5}{2}$, and we confirm that the inverse, $\frac{2}{5}$, is indeed the density of $S$ in $V(G)$.

\subsection{Cycles}

\begin{figure}[ht]
    \centering
    \begin{tabular}{c@{\hskip 5em}c}
        \begin{tabular}{c}\includegraphics[width=0.3\textwidth]{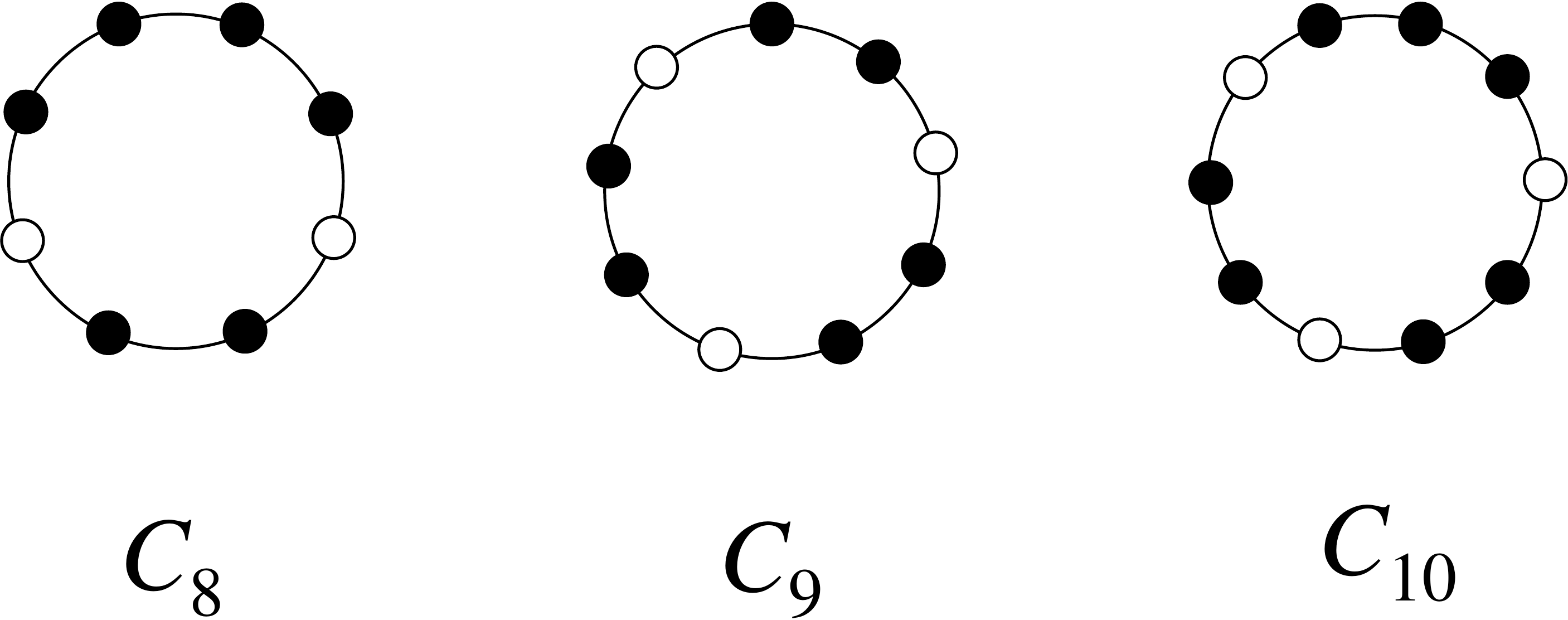}\end{tabular} &
        \begin{tabular}{c|c}
            $n \ge 5$ & $|S|$ \\ \hline
            $3k$ & $2k$ \\
            $3k + 1$ & $2k + 1$ \\
            $3k + 2$ & $2k + 2$
        \end{tabular}
    \end{tabular}
    \caption{Construction of optimal RED:LD sets, $S$, on $C_n$.}
    \label{fig:red-finite-cycle}
\end{figure}

\begin{prop}\label{theo:red-ld-finite-cycle}
RED:LD($C_n$) is $\ceil{\frac{2n}{3}}$ for $n \ge 5$, or $n$ for $n \le 4$.
\end{prop}
\begin{proof}
Let $S$ be a RED:LD set on $C_n$; by Theorem~\ref{theo:red-ld-char}, we know $S$ is 2-dominating.
As there are $n$ vertices, we require at least $2n$ dominations in total.
Because any detector can dominate at most three vertices, we require at least $\ceil{\frac{2n}{3}}$ detectors, so $|S| \ge \ceil{\frac{2n}{3}}$.
Let $V(C_n) = \{v_1,v_2,\hdots,v_n\}$ and $(v_i, v_{(i \!\!\mod n) + 1}) \in E(C_n)$.
If $n \le 4$ then RED:LD($C_n$) $= n$, as any non-detector will violate Theorem~\ref{theo:red-ld-char}; otherwise, we will assume $n \ge 5$.
Let $S = \{v_i : i \!\!\mod 3 \neq 0\}$.
Figure~\ref{fig:red-finite-cycle} shows constructions of $S$ for $C_8$, $C_9$, and $C_{10}$ along with a table of $|S|$ 
for each case of $n \!\!\mod 3$; for any case, simple algebraic manipulation shows that  $|S| = \ceil{\frac{2n}{3}}$.
Clearly, by this construction, every vertex is 2-dominated, and demonstrating that each vertex pair is distinguished follows identically to the proof of Theorem~\ref{theo:red-ld-finite-path}.
Thus, $|S|$ is an optimal RED:LD set on $C_n$.
\end{proof}

\FloatBarrier
\subsection{Ladders}

\begin{figure}[ht]
    \centering
    \includegraphics[width=0.6\textwidth]{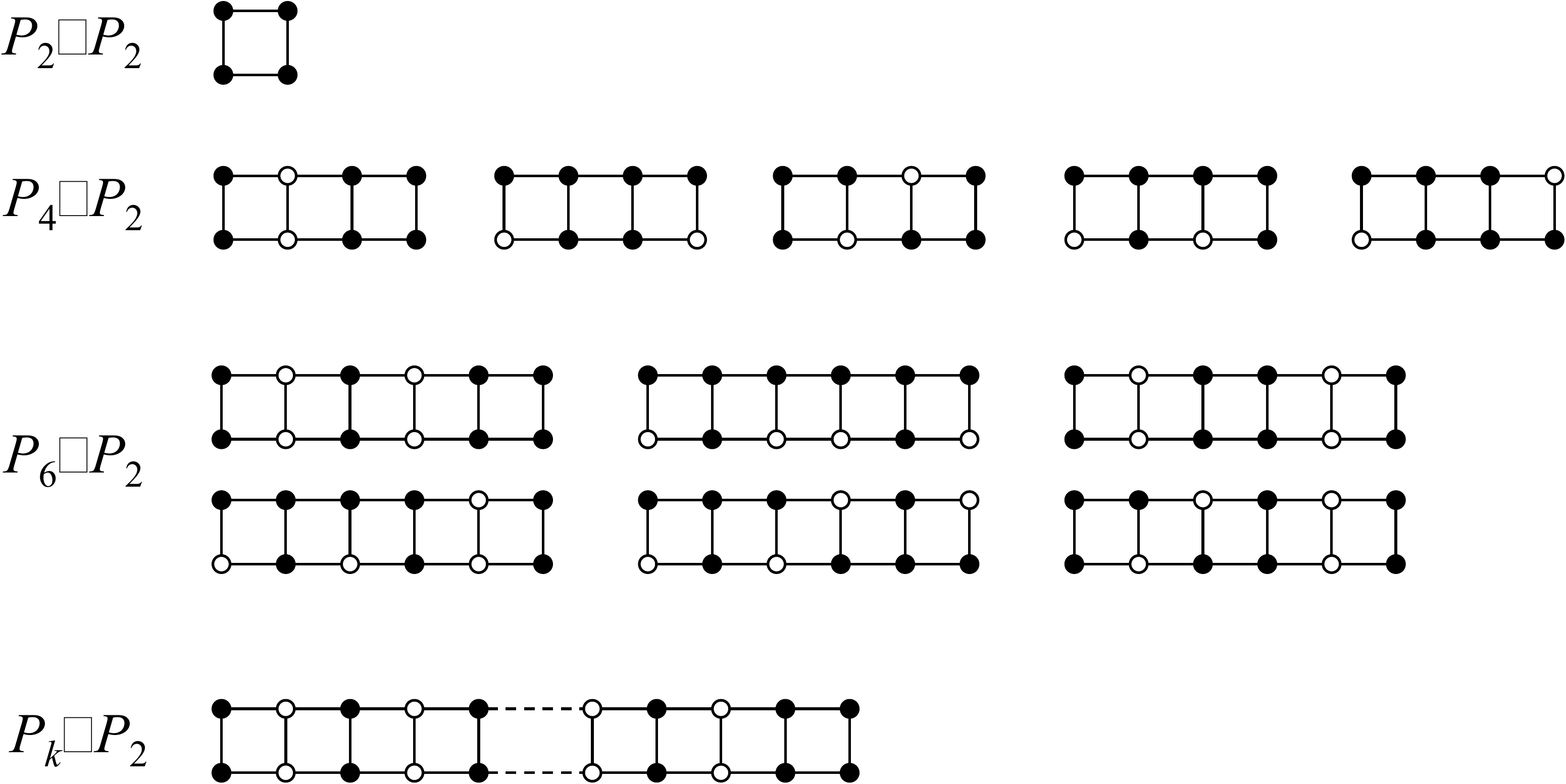}
    \caption{Optimal RED:LD sets on even ladders. For $k \le 6$, these are exhaustive.}
    \label{fig:ladders}
\end{figure}

\begin{theorem}
For a ladder graph $G = P_k \square P_2$, then RED:LD($G$) is $k+1$ if $k$ is odd, or $k+2$ otherwise.
\end{theorem}
\begin{proof}
Let $V(P_k) = \{u_1,u_2,\hdots,u_k\}$, $V(P_2) = \{w_1, w_2\}$, and $v_{i,j} \in V(G)$ denote vertex $(u_i, w_j)$.
Let $A = \{v_{i,j} : i \!\!\mod 2 = 1\} \cup \{v_{k,1}, v_{k,2}\}$; for even $k$, this construction is shown in Figure~\ref{fig:ladders} for $P_k \square P_2$.
It can be verified with Theorem~\ref{theo:red-ld-char} that $A$ is a RED:LD set, and by construction we have $|A| = k + 1$ for odd $k$ and $|A| = k + 2$ for even $k$, making these upper bounds for $\textrm{RED:LD}(G)$.

Let $S$ be a RED:LD set on $G$, which has $n = 2k$ vertices.
Every vertex must be 2-dominated, so there must be at least $2n=4k$ dominations in total.
If both vertices on one end of the graph are non-detectors, then they cannot be 2-dominated; thus, we can assume that at least one vertex on each end is a detector.
The two corner detectors will each contribute 3 dominations; the other $4k-6$ must come from the other vertices.
We know that $\Delta(G)=3$, so we require at least $\ceil{\frac{4k-6}{4}}$ additional detectors, giving a total of $|S| \ge \ceil{\frac{4k-6}{4}} + 2 = k+1$.
If $k$ is odd, this matches the upper bound, so RED:LD($G$) $= k + 1$.

For even $k$, we will use an inductive argument to show that RED:LD($G$) $= k + 2$.
If $k=2$ or $k=4$, then clearly RED:LD($G$) $= k + 2$, as shown in Figure~\ref{fig:ladders}; thus, we will assume $k \ge 6$.
Suppose that $\forall j$, $\{v_{j,1},v_{j,2}\} \cap S \neq \varnothing$.
Thus, there is at least one detector in every ``column" of $G$.
From the previously established lower bound, we know $|S| \ge k + 1$, meaning there is at least one column with two detectors.
If there is a second column with two detectors, we will have $|S| \ge k + 2$, and we are done; otherwise, we assume there is exactly one column with two detectors.
Suppose an end column has two detectors, say $\{v_{1,1},v_{1,2}\} \subseteq S$; then without loss of generality let $v_{2,1} \notin S$ and $v_{2,2} \in S$.
To distinguish $v_{1,2}$ and $v_{2,1}$, we require $v_{3,1} \in S$ and $v_{3,2} \notin S$.
To distinguish $v_{2,1}$ and $v_{3,2}$, we require $v_{4,1} \notin S$ and $v_{4,2} \in S$.
We see that $v_{3,1}$ is not 2-dominated, a contradiction.
Otherwise, both end columns have exactly one detector and without loss of generality let $v_{1,1} \in S$ and $v_{1,2} \notin S$.
To 2-dominate $v_{1,1}$ and $v_{1,2}$, we require $\{v_{2,1},v_{2,2}\} \subseteq S$.
To distinguish $v_{1,2}$ and $v_{2,1}$, we require $v_{3,1} \in S$, which implies $v_{3,2} \notin S$.
If $v_{4,1} \in S$ then we require $v_{5,2} \in S$ to 2-dominate $v_{4,2}$; we see that the two non-detectors, $v_{4,2}$ and $v_{5,1}$ cannot be distinguished, a contradiction.
Otherwise, we can assume $v_{4,2} \in S$.
To distinguish $v_{3,2}$ and $v_{4,1}$, we require $v_{5,1} \in S$; we see that $v_{4,2}$ is not 2-dominated, a contradiction.

Otherwise, we can assume $\exists j$ such that $\{v_{j,1}, v_{j,2}\} \cap S = \varnothing$.
Because each end of the graph requires at least one corner vertex, so we can assume $1 < j < k$.
To 2-dominate $v_{j,i}$, we require $\{v_{j-1,i},v_{j+1,i}\} \subseteq S$.
Consider the graph $G' = G - \{v_{j,1}, v_{j,2}\}$; because of the placement of the four detectors around $j$, it must be the case that RED:LD($G'$) $=$ RED:LD($G$).
$G'$ is now split into two pieces, $G'_1$ and $G'_2$, each of which is a ladder graph; $G'_1$ has length $j-1$, and $G'_2$ has length $k-j$.
Since $k$ is even, one of $G'_1$ or $G'_2$ is an even ladder and the other is an odd ladder.
Without loss of generality, we can assume that $G'_1$ is odd and $G'_2$ is even by flipping the graph and adjusting the value of $j$.
Because $G'_1$ is an odd ladder, RED:LD($G'_1$) $= j$; by induction, RED:LD($G'_2$) $= k-j+2$, so $\textrm{RED:LD}(G) = j + k - j + 2 = k + 2$ and we are done.
\end{proof}

\FloatBarrier
\subsection{Hypercubes}

\begin{figure}[ht]
    \centering
    \includegraphics[width=0.5\textwidth]{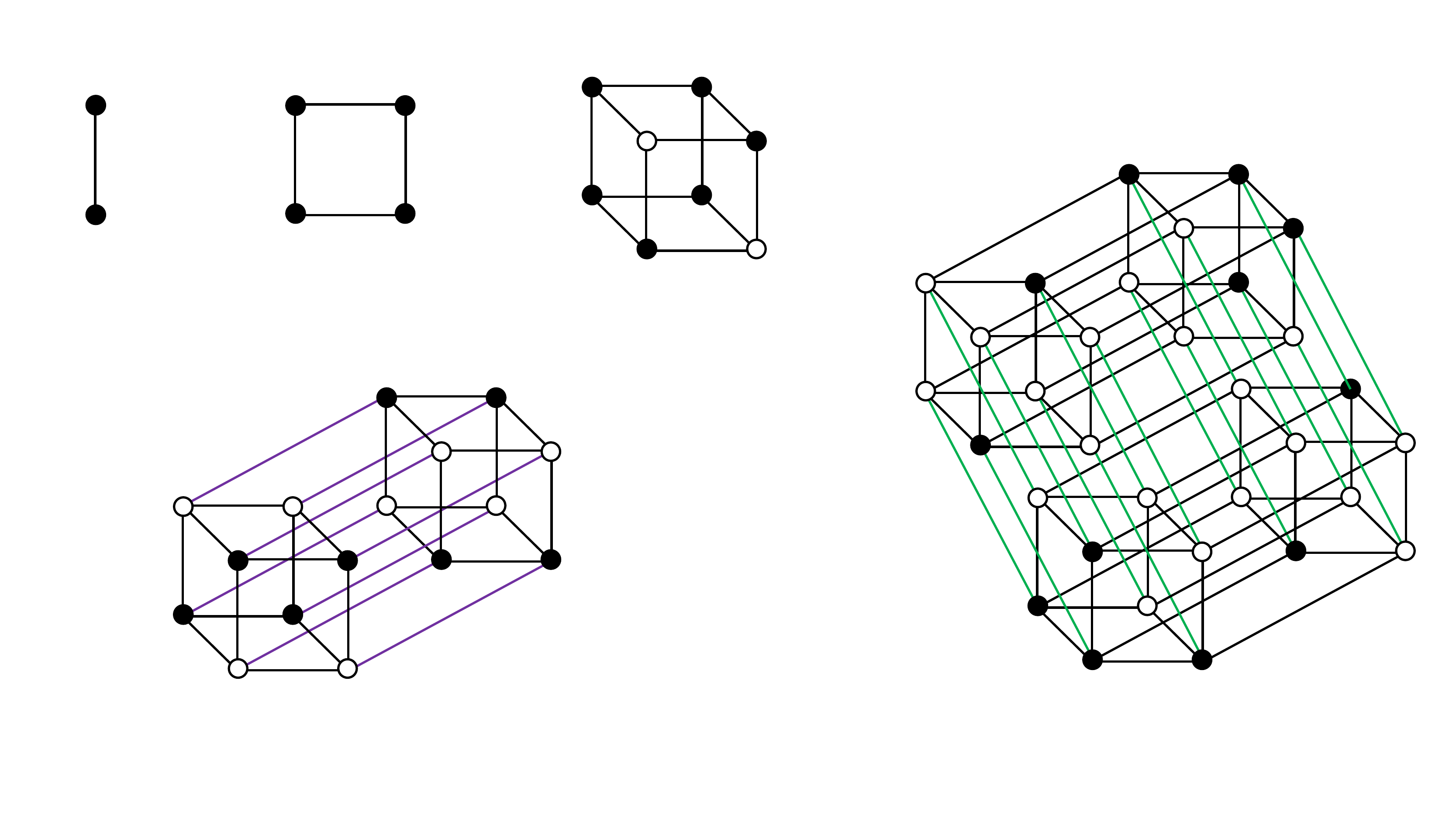}
    \caption{Optimal RED:LD sets for $Q_n$ with $n \le 5$}
    \label{fig:hypercubes}
\end{figure}

Let $Q_n = P_2^n$, where $G^n$ denotes repeated application of the $\square$ operator, be the hypercube in $n$ dimensions.
If $S$ is a RED:LD set on $Q_n$ for $n \ge 2$, then we can duplicate the vertices to produce a new RED:LD set of size $2|S|$ on $Q_{n+1} = Q_n \square P_2$; thus, $\textrm{RED:LD\%}(Q_n)$ is a non-increasing sequence in terms of $n$.
We have found that $\textrm{RED:LD\%}(Q_5) = \frac{3}{8}$, which serves as an upper bound for the minimum density of RED:LD sets in larger hypercubes.
Figure~\ref{fig:hypercubes} shows an optimal RED:LD set for each of the hypercubes on $n \le 5$ dimensions.

\FloatBarrier
\subsection{Trees}

\begin{prop}\label{prop:tree-2-dom}
If $T$ is a tree and $S \subseteq V(T)$ at least 2-dominates all vertices, then $S$ is a RED:LD set.
\end{prop}
\begin{proof}
Suppose $u \notin S$ and $v \in V(T)$; then $\exists w \in (N[u] - N[v]) \cap S$.
Thus, if $v \notin S$, then $u$ and $v$ are 2-distinguished, and if $v \in S$, then $u$ and $v$ are 1-distinguished.
By Theorem~\ref{theo:red-ld-char}, $S$ is a RED:LD set.
\end{proof}

\begin{observation}\label{theo:ceil-thing}
$a \ceil{b} \ge \ceil{ab}$ for any $a \in \mathbb{N}$ and $b \in \mathbb{R}$.
\end{observation}
\begin{proof}
If $b \in \mathbb{Z}$, then clearly $a \ceil{b} = ab = \ceil{ab}$.
Otherwise, $b = c + d$ where $c \in \mathbb{Z}$ and $d \in (0, 1)$; then $a \ceil{b} = a(c + 1) \ge \ceil{a(c + d)} = \ceil{ab}$.
\end{proof}

\begin{theorem}\label{theo:finite-tree-lower}
Let $T_n$ be a tree on $n \ge 2$ vertices; then $\textrm{RED:LD}(T_n) \ge \ceil{\frac{2n+2}{3}}$.
\end{theorem}
\begin{proof}

The proof will follow inductively; as a base case, we use $T_2$, for which the theorem holds, and we assume $n \ge 3$.
Let $S \subseteq V(T_n)$ be a RED:LD set.
If $S = V(T_n)$ then clearly $|S| \ge \ceil{\frac{2n+2}{3}}$ and we would be done; otherwise, let $v \notin S$ be some non-detector vertex.
We will break the graph into $j = |N(v)| \ge 2$ sub-trees, being the connected components of $T_n-\{v\}$; let these branches be $B_1, B_2, \hdots, B_j$ and let $n_i = |V(B_i)|$.
Because $v \notin S$, $|S \cap V(B_i)| \ge \textrm{RED:LD}(B_i) \ge \ceil{(2n_i+2)/3}$ by induction, and $|S| \ge \textrm{RED:LD}(T_n) \ge \sum_{i=1}^j{\textrm{RED:LD}(B_i)} = \sum_{i=1}^j{\ceil{(2n_i+2)/3}}$.
By applying Observation~\ref{theo:ceil-thing}, we see that $3|S| \ge 3\sum_{i=1}^j{\ceil{(2n_i+2)/3}} \ge \sum_{i=1}^j{\ceil{2n_i+2}} = 2j + 2\sum_{i=1}^j{n_i} = 2j + 2(n-1)$.
We know that $j \ge 2$ by hypothesis, so $3|S| \ge 2n + 2$; thus, $|S| \ge \frac{2n+2}{3}$; additionally, we know that $|S| \in \mathbb{N}$, so we can strengthen this to $|S| \ge \ceil{\frac{2n+2}{3}}$, completing the proof.
\end{proof}

By Observation~\ref{obs:red-ld-upper} we have $|\textrm{RED:LD}(K_{n_1,n_2})| = n_1+n_2$, and hence a star graph is a tree which requires $\textrm{RED:LD}(G) = |V(G)|$.
By Theorem~\ref{theo:red-ld-char}, we are guaranteed that a RED:LD set exists on any tree with order $n \ge 2$, and hence we have the following extremal values on RED:LD(T) for a tree T of order $n$.

\begin{corollary}
Let $T_n$ be a tree of order $n \ge 2$; then $\ceil{\frac{2n+2}{3}} \leq \textrm{RED:LD}(T_n) \leq n$.
\end{corollary}

\subsubsection{Extremal trees with $\textrm{RED:LD}(T_n) = n$}

\begin{theorem}\label{theo:redld-max-tree-leaf-support}
A tree, $T$, has $\textrm{RED:LD}(T)=n$ if and only if every vertex is a leaf or support vertex.
\end{theorem}
\begin{proof}
Because a RED:LD set, $S$, must at least 2-dominate all vertices, we see that any leaf vertex, $v \in V(T)$ requires $N[v] \subseteq S$; thus, if every vertex is a leaf or support vertex, then $\textrm{RED:LD}(T)=n$.
For the converse, suppose $\exists v \in V(T)$ which is not a leaf or support vertex; let $S = V(T)-\{v\}$.
Because $v$ is not a leaf vertex, $deg(v) \ge 2$, meaning $v$ is at least 2-dominated by $S$.
Every vertex outside of $N[v]$ will be at least 2-dominated because $T$ is connected and $v$ is the only non-detector.
Any vertex $u \in N(v)$ is at least 1-dominated by itself; if it is only 1-dominated, then $deg(u)=1$, which contradicts that $v$ is not a support vertex.
Thus, $S$ causes all vertices to be at least 2-dominated, and Proposition~\ref{prop:tree-2-dom} yields that $S = V(T)-\{v\}$ is a RED:LD set.
Therefore, if there is a non-leaf non-support vertex, then $\textrm{RED:LD}(T) < n$, completing the proof.
\end{proof}

Let $\mathscr{T}_{max}$ denote the family of all trees of order $n$ with $\textrm{RED:LD}(T_n) = n$.
We will now show how we can generate the entire set $\mathscr{T}_{max}$.

\begin{theorem}\label{theo:redld-max-add-1}
Let $T \in \mathscr{T}_{max}$ with $u \in V(T)$, and let $v$ be a new vertex.
Then $T' = (V(T) \cup \{v\}, E(T) \cup \{vu\}) \in \mathscr{T}_{max}$ if and only if $u$ is a support vertex or a leaf where its support vertex has at least two leaves.
\end{theorem}
\begin{proof}
From Theorem~\ref{theo:redld-max-tree-leaf-support}, we know that every vertex is either a support or leaf vertex.
Clearly, if $u$ is a support vertex, then adding $v \in N(u)$ will still result in every vertex being either a leaf or support, so Theorem~\ref{theo:redld-max-tree-leaf-support} has that $T' \in \mathscr{T}_{max}$.
Suppose $u$ is a leaf where its support vertex, $w \in N(u)$, has at least two leaves.
By adding $v \in N(u)$, $w$ remains a support vertex due to its remaining leaf, and $u$ switches from being a leaf to being a support vertex; every vertex remains a leaf or support vertex, so $T' \in \mathscr{T}_{max}$.
For the converse, suppose that $u$ is a leaf which is the only leaf of its adjacent support vertex, $w$.
Then adding $v \in N(u)$ causes $w$ to no longer be a leaf or support vertex.
Theorem~\ref{theo:redld-max-tree-leaf-support} gives that $T' \notin \mathscr{T}_{max}$, completing the proof.
\end{proof}

\begin{theorem}\label{theo:redld-max-remove-1}
Let $T \in \mathscr{T}_{max}$ with $n \ge 3$, and let $v \in V(T)$.
Then $T' = T - \{v\} \in \mathscr{T}_{max}$ if and only if $v$ is a leaf with at least one other sibling leaf or $v$ is a leaf where its support vertex has degree 2.
\end{theorem}
\begin{proof}
From Theorem~\ref{theo:redld-max-tree-leaf-support}, we know that every vertex is either a support or leaf vertex.
Clearly, if $v$ is a leaf with at least one other sibling leaf, then removing $v$ still results in every vertex being a leaf or support vertex, so $T' \in \mathscr{T}_{max}$.
Suppose $v$ is a leaf where its support vertex, $u \in N(v)$ has $deg(u) = 2$.
By removing $v$, vertex $u$ goes from being a support vertex to a leaf vertex; every vertex remains a leaf or support, so $T' \in \mathscr{T}_{max}$.
For the converse, suppose that $v$ is either a support vertex, or a leaf vertex which is the only leaf of its support vertex $u \in N(v)$ and $deg(u) \ge 3$.
Clearly, if $v$ is a support vertex, then $T'$ is not a tree due to not being connected, so we assume the second possibility.
We see that removing $v$ causes $u$ to no longer be a support vertex, and $u$ is not a leaf because $deg(u) \ge 3$.
Thus, $u$ is neither a support nor a leaf vertex, so $T' \notin \mathscr{T}_{max}$, completing the proof.
\end{proof}

\begin{figure}[ht]
    \centering
    \includegraphics[width=0.99\textwidth]{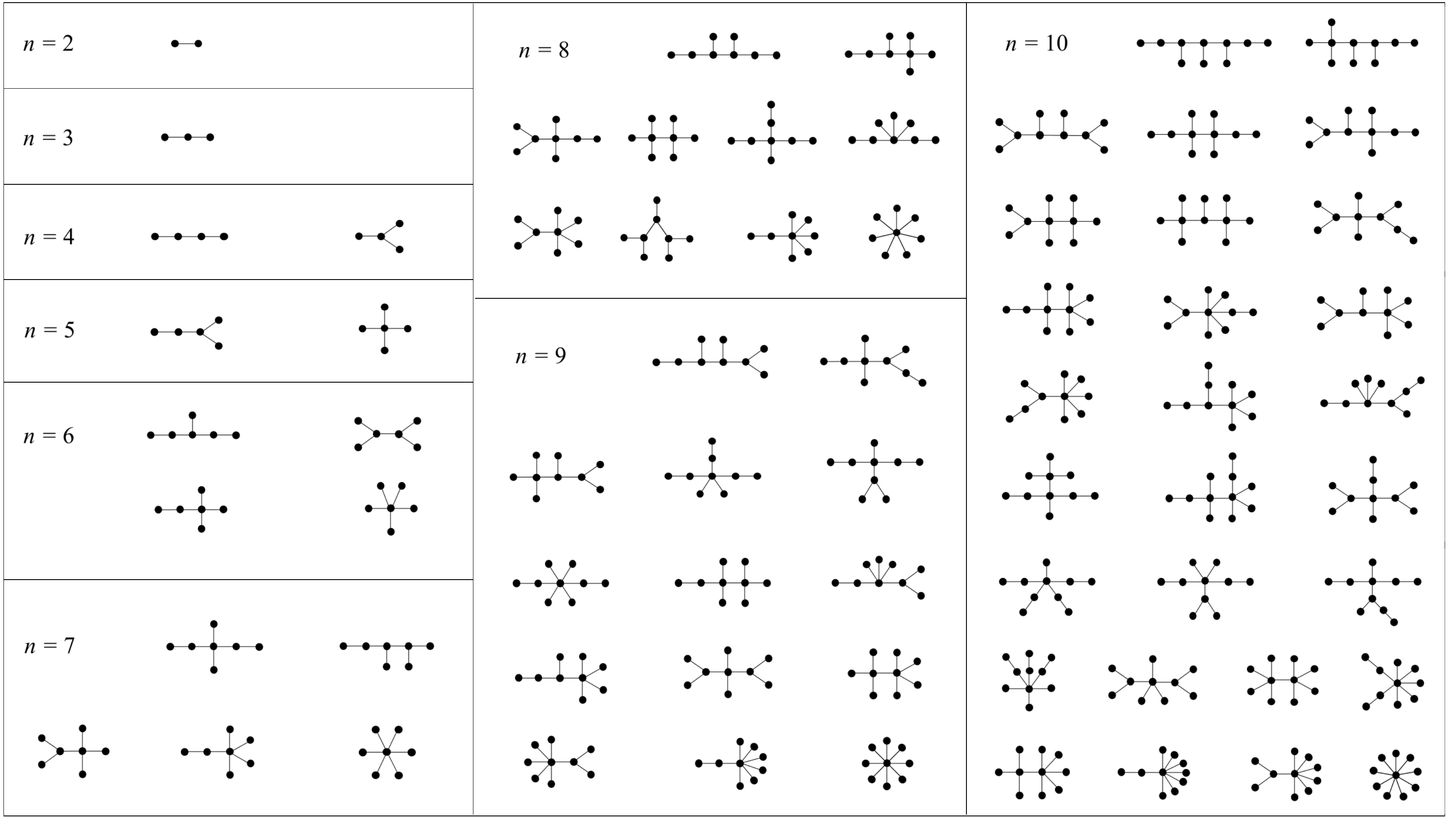}
    \caption{All trees in $\mathscr{T}_{max}$ of order $n \le 10$}
    \label{fig:redld-max-trees-all}
\end{figure}

\begin{theorem}\label{theo:redld-tree-max-complete}
If $T \in \mathscr{T}_{max}$ has $n \ge 3$, then there is a leaf vertex $v \in V(T)$ such that $T-v \in \mathscr{T}_{max}$.
\end{theorem}
\begin{proof}
By Theorem~\ref{theo:redld-max-tree-leaf-support}, we know that every vertex in $T$ is either a leaf or support vertex, and because $n \ge 3$, we know that there is at least one support vertex.
Let $T'$ be the graph generated by the set of support vertices, let $u \in V(T')$ be a leaf vertex in $T'$, and let $v \in N(u)$ by one of its leaves in the original graph, $T$.
Then $v$ satisfies the requirements of Theorem~\ref{theo:redld-max-remove-1}, so $T-v \in \mathscr{T}_{max}$, completing the proof.
\end{proof}

By repeatedly applying Theorem~\ref{theo:redld-tree-max-complete} to an arbitrary tree $T \in \mathscr{T}_{max}$, we see that we will eventually hit the $P_2$ base case.
By performing the vertex removal steps in reverse, we see that $T$ can be constructed from the base case tree, $P_2$, by adding vertices one at a time, with every intermediate tree being likewise in $\mathscr{T}_{max}$.
Thus, the construction process given by Theorem~\ref{theo:redld-max-add-1} constructs the entire family $\mathscr{T}_{max}$.
Figure~\ref{fig:redld-max-trees-all} shows all trees in $\mathscr{T}_{max}$ on $n \le 10$ vertices.



\FloatBarrier
\subsubsection{Extremal trees with $\textrm{RED:LD}(T_n) = \ceil{\frac{2n+2}{3}}$}

\begin{figure}[ht]
    \centering
    \begin{tabular}{|c|c|} \hline
        \includegraphics[width=0.4\textwidth]{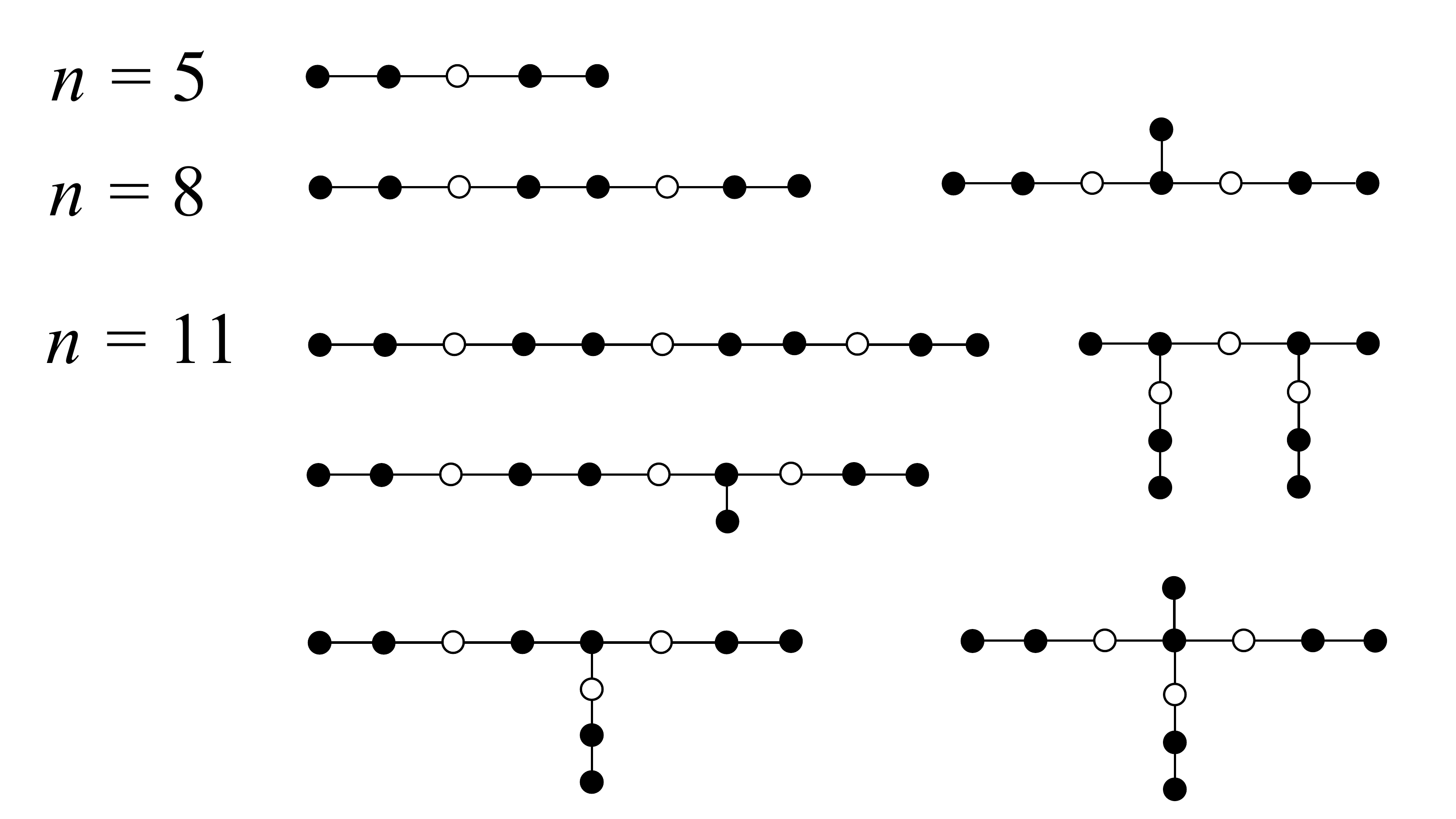} &  
        \includegraphics[width=0.4\textwidth]{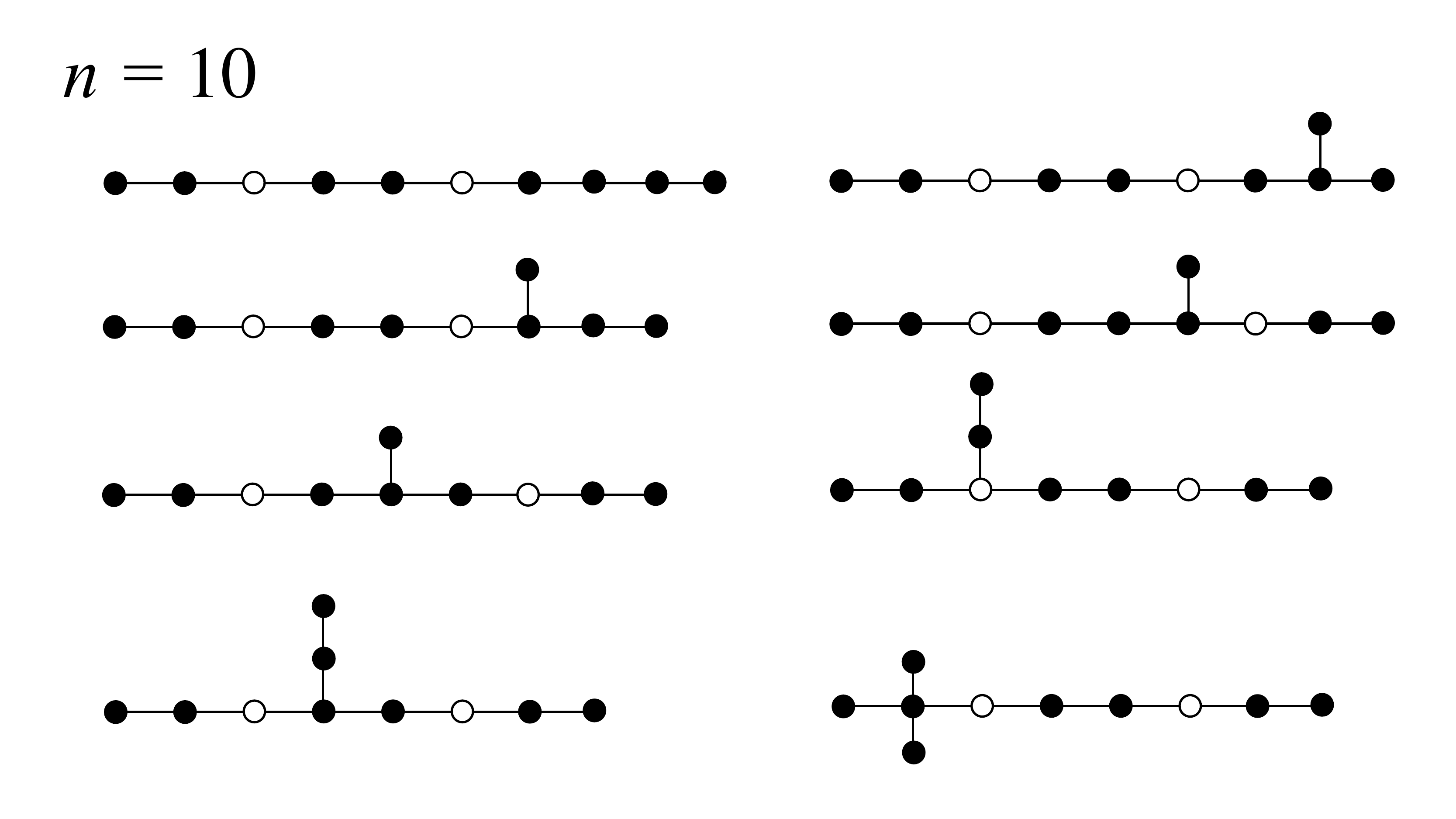} \\
        Trees in $\mathscr{T}_{min}^2$ for $n=5,8,11$ & \\\cline{1-1}
        \includegraphics[width=0.4\textwidth]{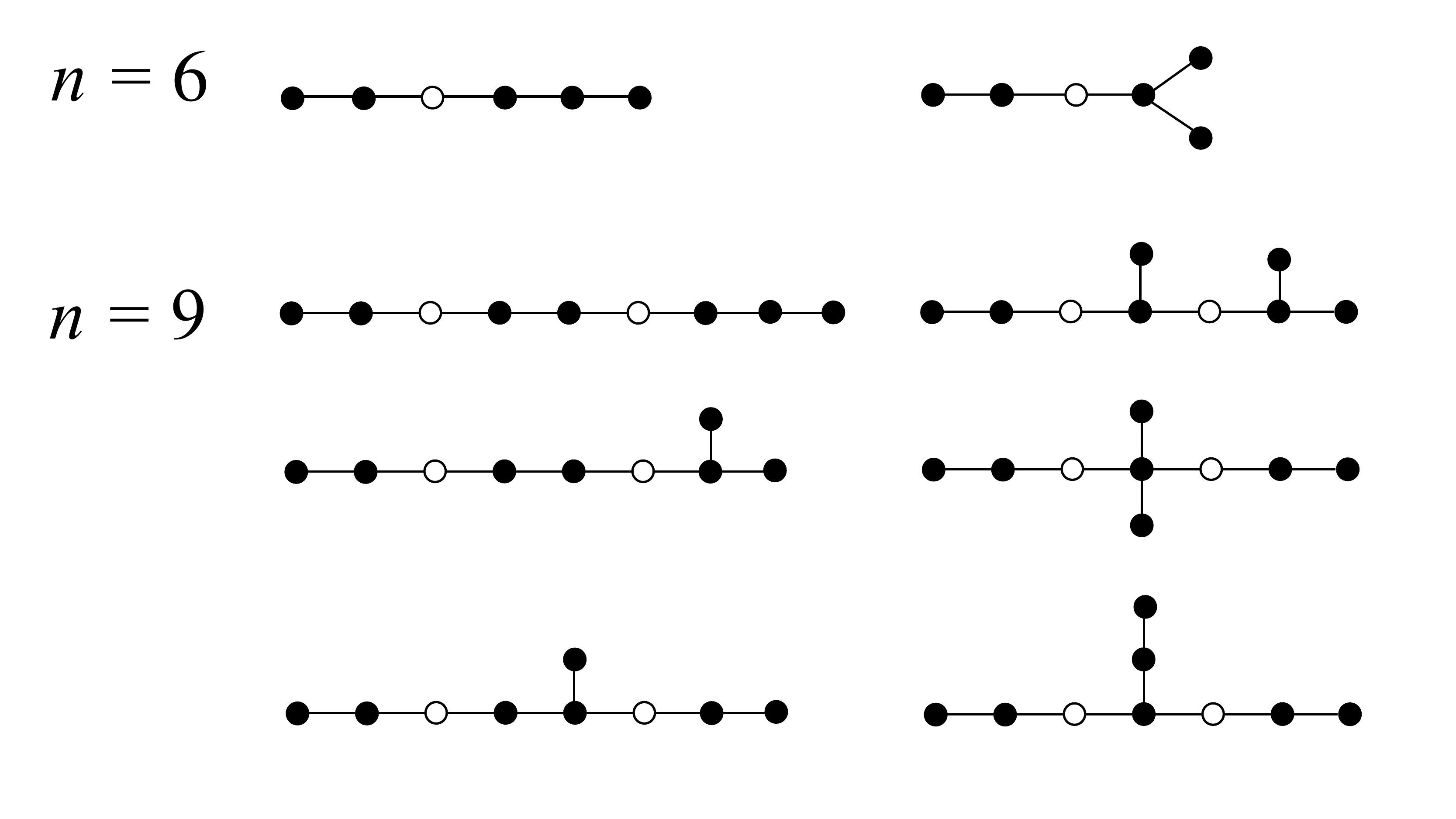} &
        \includegraphics[width=0.4\textwidth]{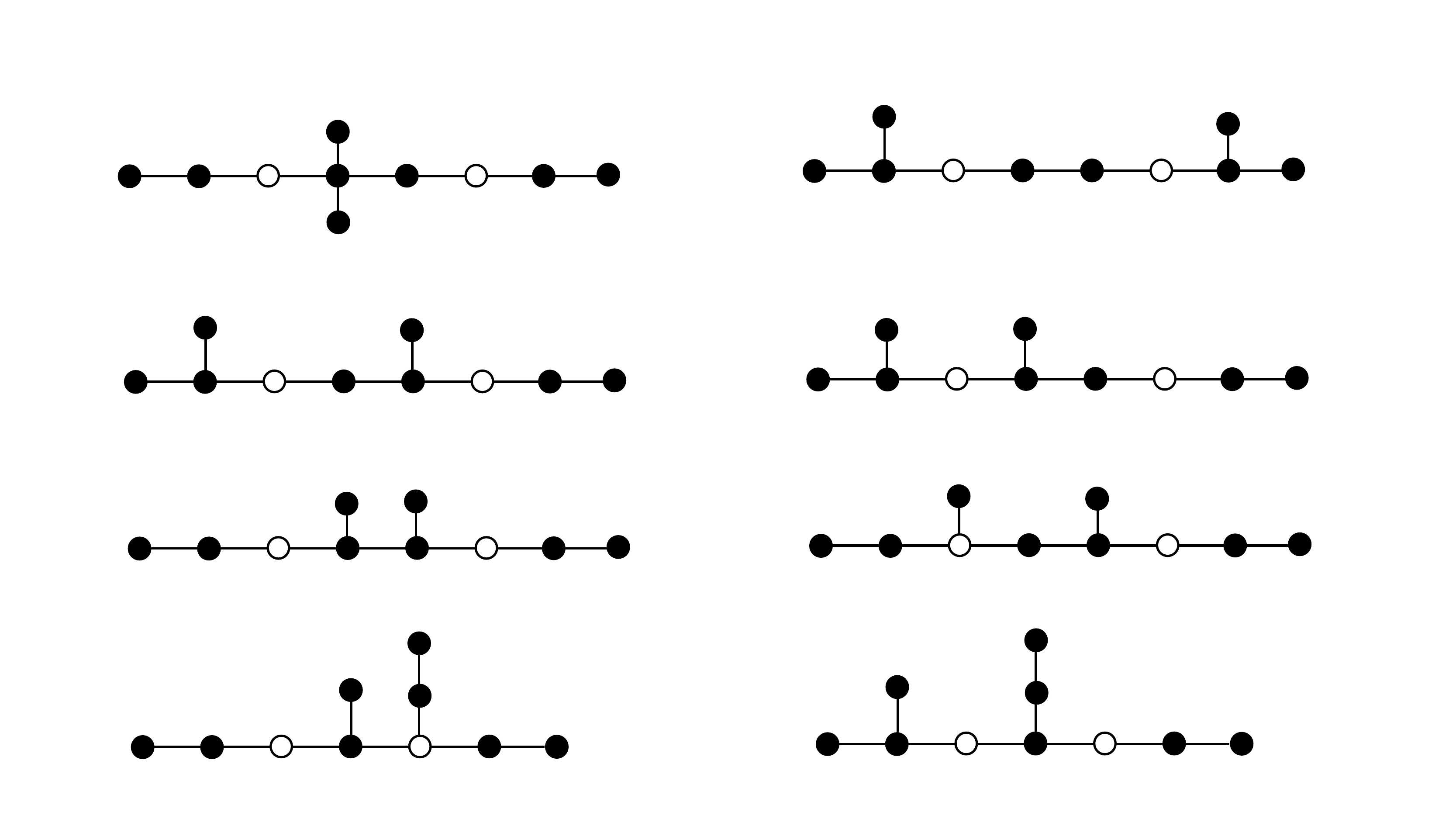} \\
        Trees in $\mathscr{T}_{min}^0$ for $n=6,9$ & \\\cline{1-1}
    \end{tabular}
    \begin{tabular}{|cc|}
        \includegraphics[width=0.4\textwidth]{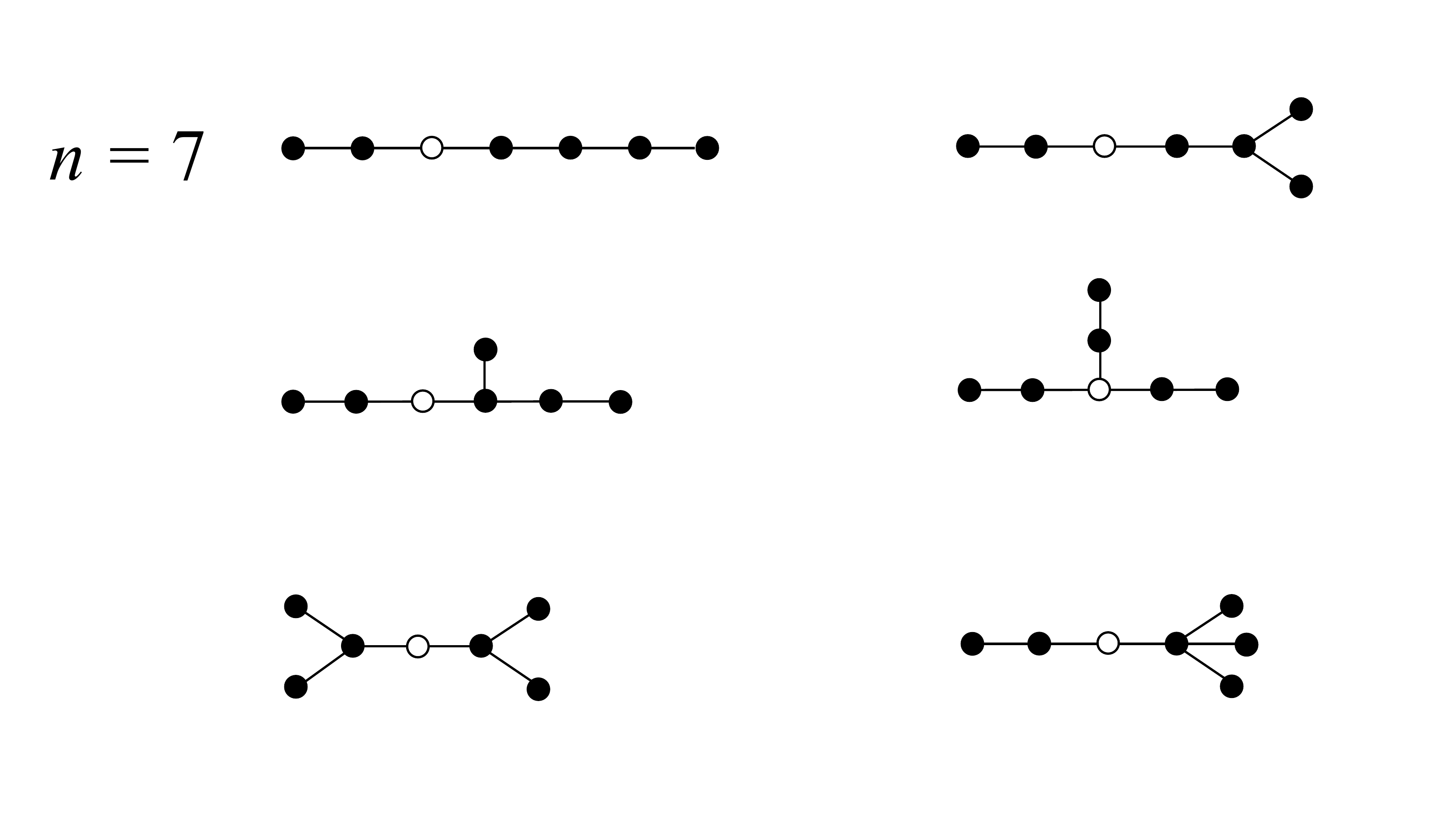} &  \includegraphics[width=0.4\textwidth]{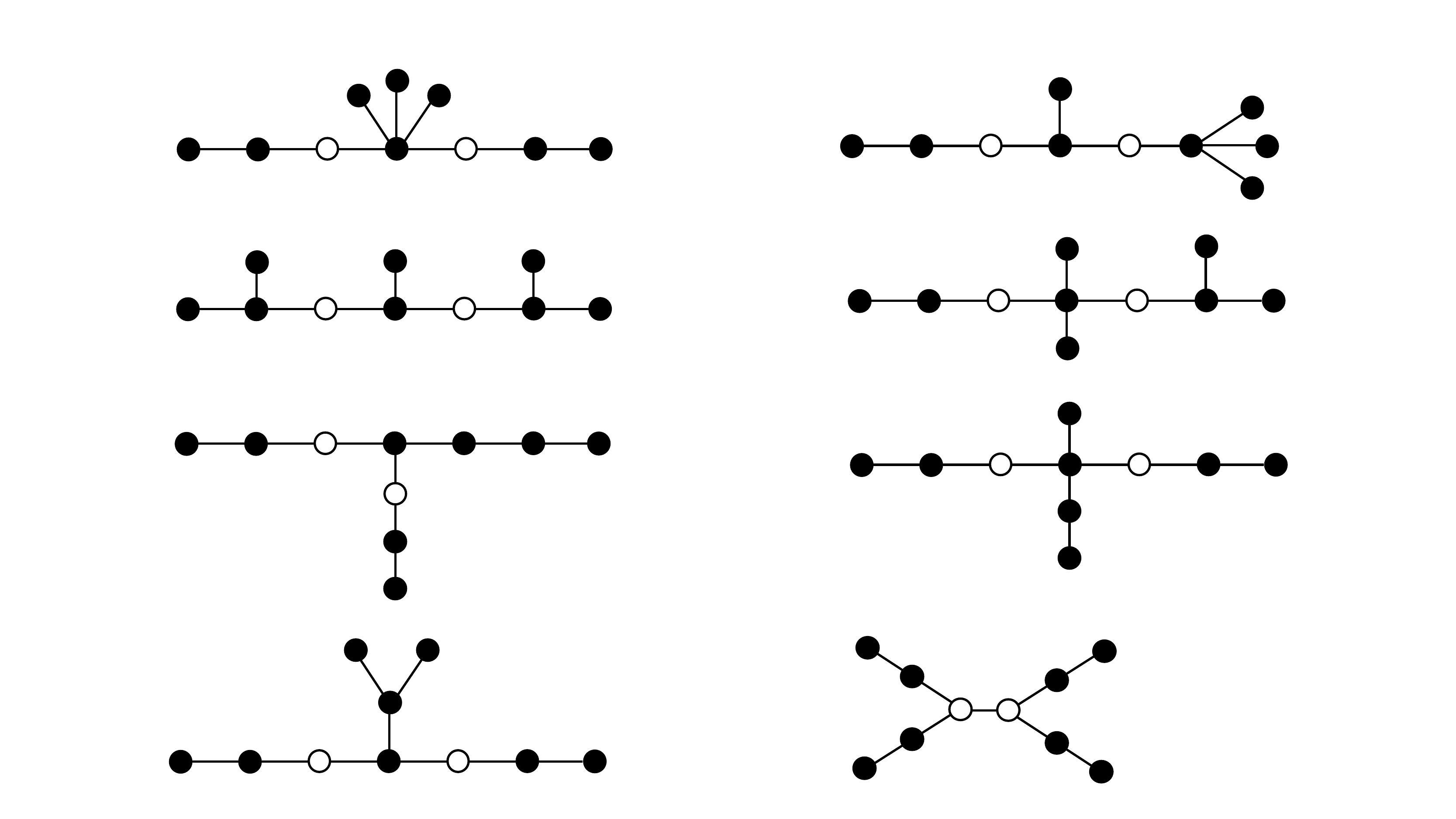} \\
        \multicolumn{2}{|c|}{Trees in $\mathscr{T}_{min}^1$ for $n=7,10$} \\\hline
    \end{tabular}
    \caption{All trees in $\mathscr{T}_{min}$ of order $5 \le n \le 11$}
    \label{fig:redld-min-trees-all}
\end{figure}

\begin{theorem}\label{theo:red-ld-finite-path}
RED:LD($P_n$) $= \ceil{\frac{2n+2}{3}}$.
\end{theorem}
\cbeginproof
The lower bound is proven by Theorem~\ref{theo:finite-tree-lower}.
Let $V(P_n) = \{v_1,v_2,\hdots,v_n\}$ and $E(P_n) = \{(v_i,v_j) : |i-j| = 1\}$.
Let $S = \{v_i : i \!\!\mod 3 \neq 0\} \cup \{v_{n-1},v_{n}\}$.
Figure~\ref{fig:redld-min-trees-all} includes constructions of $S$ for $P_8$, $P_9$, and $P_{10}$, and Table~\ref{tab:redld-min-tree} gives table of values of $|S|$ for each case of $n \!\!\mod 3$; for any case, simple algebraic manipulation shows that $|S| = \ceil{\frac{2n+2}{3}}$.
From this construction, it is clear that every vertex is 2-dominated; thus, by Proposition~\ref{prop:tree-2-dom}, $S$ is an optimal RED:LD set, completing the proof.
\cendproof

\begin{table}[ht]
    \centering
    \begin{tabular}{c|c||c|c}
        $n$ & $|S|$ & $j = n-|S|$ & $n$\\ \hline
        $3k$ & $2k + 1$ & $k-1$ & $3j+3$ \\
        $3k + 1$ & $2k + 2$ & $k-1$ & $3j+4$ \\
        $3k + 2$ & $2k + 2$ & $k$ & $3j + 2$
    \end{tabular}
    \caption{Lower bound RED:LD values for $T_n$}
    \label{tab:redld-min-tree}
\end{table}


\begin{prop}\label{theo:red-ld-inf-path}
RED:LD($P_\infty$) $= \frac{2}{3}$.
\end{prop}
\begin{proof}
Let $x \in S \subseteq V(P_\infty)$.
From construction, $|N[x]| = 3$, and by Theorem~\ref{theo:red-ld-char} we know every vertex must be at least 2-dominated.
Thus, $sh(x) \le 3 \times \frac{1}{2}$, giving a lower bound density of $\frac{2}{3}$.
Where $V(P_\infty) = \{v_k : k \in \mathbb{Z}\}$, let $S = \{v_k : k \!\!\mod 3 \neq 0\}$; then the density of $S$ in $V(G)$ is $\frac{2}{3}$.
By this construction $S$ is clearly 2-dominating, so Proposition~\ref{prop:tree-2-dom} has that $S$ is an optimal RED:LD set.
\end{proof}

\begin{wrapfigure}[9]{r}{0.35\textwidth}
    \centering
    \includegraphics[width=0.3\textwidth]{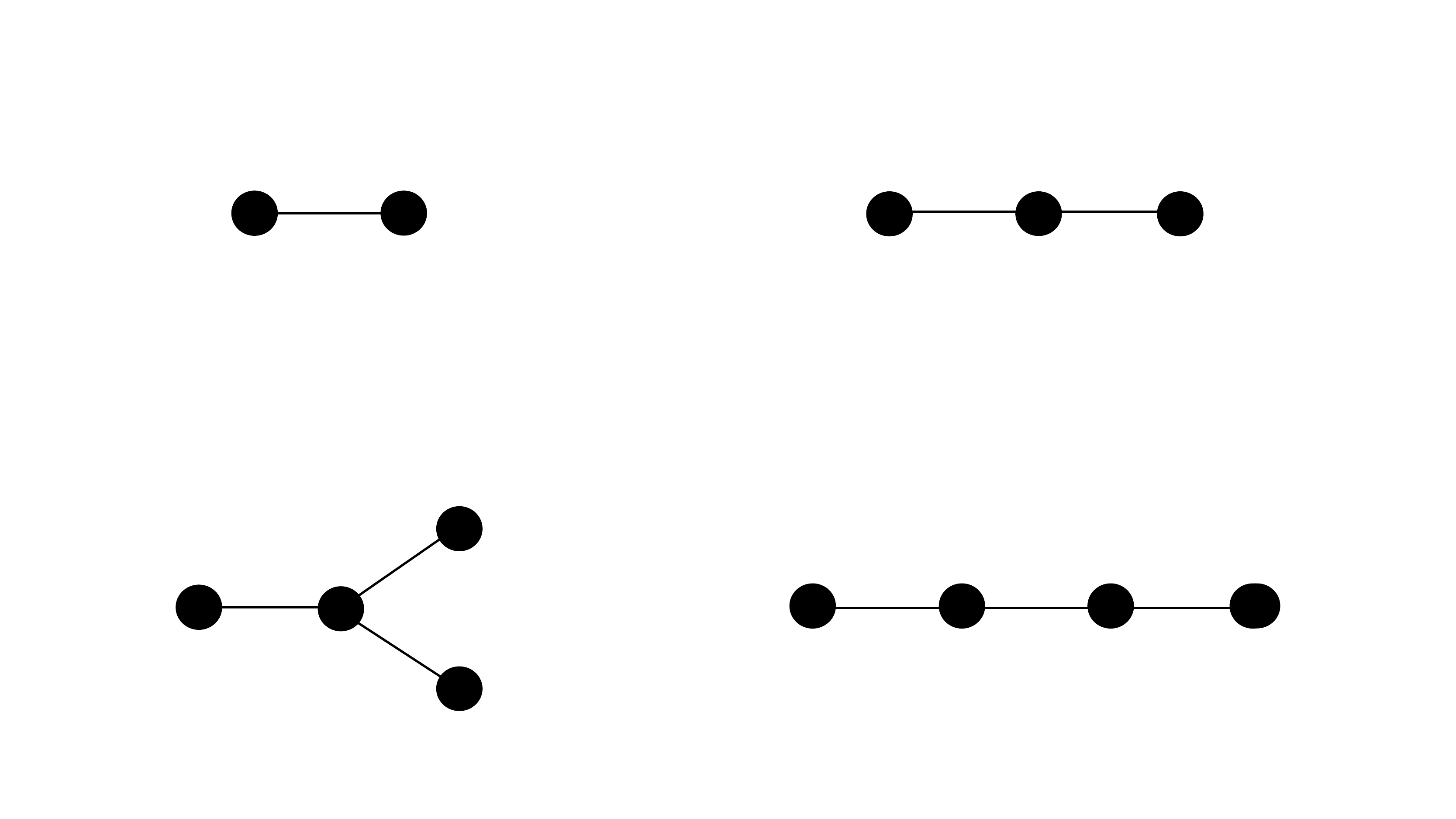}
    \caption{Elements of $\mathscr{T}_{min}$ for $n \le 4$}
    \label{fig:redld-extreme-low-tree}
\end{wrapfigure}

From Theorems \ref{theo:finite-tree-lower} and \ref{theo:red-ld-finite-path} we see that among all trees, paths have the lowest value of $\textrm{RED:LD}(G)$.
These values are broken down for all $n$ mod 3 in Table~\ref{tab:redld-min-tree}.

In what follows, we let $\mathscr{T}_{min}$ denote the family of all trees of order $n$ with $\textrm{RED:LD}(T_n) = \ceil{\frac{2n+2}{3}}$ and we let $S$ be an optimal RED:LD set of $T_n$.
Further, we let $\mathscr{T}_{min}^i \subseteq \mathscr{T}_{min}$ denote the trees of order $3k+i$.
The trees on up to four vertices in $\mathscr{T}_{min}$ are given in Figure~\ref{fig:redld-extreme-low-tree}.
We will now show rules that can be used to generate the entire set $\mathscr{T}_{min}$.

\vspace{1em}
\begin{observation}\label{obs:max-n-for-k-nondetector}
Let $S$ be a RED:LD set for a tree, $T$, and $j = |V(T)-S|$; then, $|V(T)| \ge 3j + 2$.
\end{observation}
\begin{proof}
Columns 3 and 4 of Table~\ref{tab:redld-min-tree} give expressions for $n$ in terms of $j$ for extremal trees of any order $n$.
Column 4 shows the smallest value for $n$ is $3j+2$, so if $T$ is extremal then we are done.
Otherwise $T$ is not extremal, in which case there must be even more (detector) vertices, completing the proof.
\end{proof}

\begin{observation}
If $T \in \mathscr{T}_{min}$ has an optimal RED:LD set with $j$ non-detectors, then $|V(T)| \le 3j + 4$.
\end{observation}
\begin{proof}
Suppose $|V(T)| \ge 3j + 5$; then it does not match any row of Table~\ref{tab:redld-min-tree}, and hence cannot be extremal.
\end{proof}

\begin{lemma}\label{lem:type-2-extremal-1}
Let $T_1, T_2 \in \mathscr{T}_{min}^2$ be trees with minimum RED:LD sets $S_1,S_2$, respectively.
Let $T$ be a tree obtained by adding a vertex, $v$, and edges $vw_1,vw_2$ where $w_1 \in S_1$ and $w_2 \in S_2$.
Then, $T \in \mathscr{T}_{min}^2$.
\end{lemma}
\begin{proof}
Because $T_1, T_2 \in \mathscr{T}_{min}^2$, we let $|V(T_1)| = 3k_1+2$ and $|V(T_2)| = 3k_2+2$ for some $k_1, k_2 \in \mathbb{N}_0$, and from Table~\ref{tab:redld-min-tree} we have $|S_1| = 2k_1 + 2$ and $|S_2| = 2k_2 + 2$.
In the combined tree, $T$, we have $|V(T)|= |V(T_1)| + |V(T_2)| + 1 = (3k_1+2) + (3k_2+2) + 1 = 3(k_1 +k_2 +1) +2$.
Let $S = S_1 \cup S_2$; then $|S| = |S_1| + |S_2| = (2k_1 +2) + (2k_2 +2) = 2(k_1+k_2+1) +2$.
Every vertex in $T$ is 2-dominated in $S$, so Proposition~\ref{prop:tree-2-dom} yields that $S$ is a RED:LD set on $T$; thus, $T \in \mathscr{T}_{min}^2$.
\end{proof}

\begin{lemma}\label{lem:type-2-extremal-2}
Let $S$ be an optimal RED:LD set for a tree $T \in \mathscr{T}_{min}^2$ with order at least 5. Then, every vertex $v \notin S$ has degree 2 and the two connected components, $T_1$ and $T_2$, of $T-v$ are in $\mathscr{T}_{min}^2$, and $S_i = S \cap V(T_i)$ is optimal on $T_i$.
\end{lemma}
\begin{proof}
Since $T$ is in $\mathscr{T}_{min}^2$, let $|V(T)| = 3k+2$ with $|S| = 2k+2$, meaning there are $j = k$ non-detectors.
Let $v \notin S$ be an arbitrary non-detector with $deg(v) = p$ and let $T_1, \hdots, T_p$ be the connected components of $T-v$ with each $T_i$ having $j_i$ non-detectors, where $j_1+\cdots+j_p=j-1$.
Let $S_i = S \cap V(T_i)$ for $1 \leq i \leq p$. 
Since $v \notin S$, $S_i$ must be a RED:LD set for $T_i$.
By Observation~\ref{obs:max-n-for-k-nondetector}, we have $|V(T_i)| \ge 3j_i+2$. Then,
\begin{align*}
|V(T)| &=  |V(T_1)| +|V(T_2)| +\cdots+ |V(T_p)| + 1 \\
&\ge  (3j_1+2) + (3j_2+2) + \cdots + (3j_p+2) + 1 \\
&= 3(j_1+j_2 +\cdots+j_p) + 2p +1 \\
&= 3(j - 1) + 2p +1
\end{align*}
Because $v$ is a non-detector, $deg(v) \ge 2$.
If $deg(v) \ge 3$, then  $|V(T)| \ge 3(j-1) + 2 \times 3 + 1 = 3k+4$, a contradiction.
Otherwise $deg(v) = 2$, and because $T$ is assumed to be in $\mathscr{T}_{min}^2$, we know $|V(T)| = 3j + 2$.
Let $|V(T_i)| = 3j_i + \alpha_i$; Observation~\ref{obs:max-n-for-k-nondetector} gives us that $\alpha_i \ge 2$.
We know that $|V(T_1)| + |V(T_2)| + 1 = |V(T)|$, which means $(3j_1 + \alpha_1) + (3j_2 + \alpha_2) + 1 = 3(j_1 + j_2) + \alpha_1 + \alpha_2 + 1 = 3(j-1) + \alpha_1 + \alpha_2 + 1 = 3j + 2$.
Thus, we see that $\alpha_1 + \alpha_2 = 4$, so $\alpha_1=\alpha_2=2$; from Table~\ref{tab:redld-min-tree}, we see these correspond to $T_1,T_2 \in \mathscr{T}_{min}^2$, and $S_1,S_2$ are optimal, completing the proof.
\end{proof}

\begin{lemma}\label{lem:type-0-extremal-1}
Let $T_1 \in \mathscr{T}_{min}^0$ and $T_2 \in \mathscr{T}_{min}^2$ be trees with minimum RED:LD sets $S_1,S_2$, respectively.
Let $T$ be a tree obtained by adding a vertex, $v$, and edges $vw_1,vw_2$ where $w_1 \in S_1$ and $w_2 \in S_2$.
Then, $T \in \mathscr{T}_{min}^0$.
\end{lemma}
\begin{proof}
Because $T_1 \in \mathscr{T}_{min}^0$ and $T_2 \in \mathscr{T}_{min}^2$, we let $|V(T_1)| = 3k_1$ and $|V(T_2)| = 3k_2+2$ for some $k_1, k_2 \in \mathbb{N}_0$, and from Table~\ref{tab:redld-min-tree} we have $|S_1| = 2k_1 + 1$ and $|S_2| = 2k_2 + 2$.
In the combined tree, $T$, we have $|V(T)|= |V(T_1)| + |V(T_2)| + 1 = (3k_1) + (3k_2+2) + 1 = 3(k_1 +k_2 +1)$.
Let $S = S_1 \cup S_2$; then $|S| = |S_1| + |S_2| = (2k_1 +1) + (2k_2 +2) = 2(k_1+k_2+1) +1$.
Every vertex in $T$ is 2-dominated in $S$, so Proposition~\ref{prop:tree-2-dom} yields that $S$ is a RED:LD set on $T$; thus, $T \in \mathscr{T}_{min}^0$.
\end{proof}

\begin{lemma}\label{lem:type-0-extremal-2}
Let $S$ be an optimal RED:LD set for a tree $T \in \mathscr{T}_{min}^0$ with order at least 6. Then, every vertex $v \notin S$ has degree 2 and the two connected components, $T_1$ and $T_2$ of $T-v$ satisfy $T_1 \in \mathscr{T}_{min}^0$ and $T_2 \in \mathscr{T}_{min}^2$, and $S_i = S \cap V(T_i)$ is optimal on $T_i$.
\end{lemma}
\begin{proof}
Since $T$ is in $\mathscr{T}_{min}^0$, let $|V(T)| = 3k$ with $|S| = 2k+1$, meaning there are $j = k\!-\!1$ non-detectors.
Let $v \notin S$ with $deg(v) = p$ and define $T_1,\hdots,T_p$, $j_1,\hdots,j_p$, and $S_1,\hdots,S_p$ similarly to Lemma~\ref{lem:type-2-extremal-2}.
We again find that $|V(T)| \ge 3(j - 1) + 2p +1$.
If $deg(v) \ge 3$, then  $|V(T)| \ge 3(j-1) + 2 \times 3 + 1 = 3k+1$, a contradiction.
Otherwise $deg(v) = 2$, and because $T$ is assumed to be in $\mathscr{T}_{min}^0$, we know $|V(T)| = 3j + 3$ (see Table~\ref{tab:redld-min-tree}).
Let $|V(T_i)| = 3j_i + \alpha_i$; Observation~\ref{obs:max-n-for-k-nondetector} gives us that $\alpha_i \ge 2$.
We know that $|V(T_1)| + |V(T_2)| + 1 = |V(T)|$, which means $(3j_1 + \alpha_1) + (3j_2 + \alpha_2) + 1 = 3(j_1 + j_2) + \alpha_1 + \alpha_2 + 1 = 3(j-1) + \alpha_1 + \alpha_2 + 1 = 3j + 3$.
Thus, we see that $\alpha_1 + \alpha_2 = 5$, so without loss of generality $\alpha_1=2$ and $\alpha_2=3$; from Table~\ref{tab:redld-min-tree}, we see these correspond to $T_1 \in \mathscr{T}_{min}^0$ and $T_2 \in \mathscr{T}_{min}^2$, and $S_1,S_2$ are optimal on $T_1,T_2$, completing the proof.
\end{proof}

\begin{lemma}\label{lem:type-1-extremal-1}
Let $\beta = (\beta_1, \hdots, \beta_p)$ with $\beta \in \{(0, 0), (1, 2), (2, 2, 2)\}$, let $T_1, \hdots, T_p$ be trees such that $T_i \in \mathscr{T}_{min}^{\beta_i}$, and let $S_i$ be a minimum RED:LD set on $T_i$.
Let $S = \cup_i{S_i}$, and let $w$ be a new vertex with $N(w) = \{x \in V(T_i)\}$ such that $|N(w) \cap S| \ge 2$.
Then the combined tree $T$ with $V(G) = \cup_i{V(T_i)} \cup \{w\}$ is in $\mathscr{T}_{min}^{1}$.
\end{lemma}
\begin{proof}
Let $|V(T_i)| = 3k_i + \beta_i$; then the combined tree, $T$, has $|V(T)| = |V(T_1)| + \cdots + |V(T_p)| + 1 = (3k_1 + \beta_1) + \cdots + (3k_p + \beta_p) + 1$.
Suppose $\beta = (0, 0)$; then $|V(T)| = 3(k_1 + k_2) + 1$ and $|S| = (2k_1 + 1) + (2k_2 + 1) = 2(k_1 + k_2) + 2$.
Suppose $\beta = (1, 2)$; then $|V(T)| = 3(k_1 + k_2 + 1) + 1$ and $|S| = (2k_1 + 2) + (2k_2 + 2) = 2(k_1 + k_2 + 1) + 2$.
Lastly, suppose $\beta = (2,2,2)$; then $|V(T)| = 3(k_1 + k_2 + k_3 + 2) + 1$ and $|S| = (2k_1 + 2) + (2k_2 + 2) + (2k_3 + 2) = 2(k_1 + k_2 + k_3 + 2) + 2$.
By hypothesis, we see that $S$ 2-dominates every vertex in $T$, so Proposition~\ref{prop:tree-2-dom} yields that $S$ is a RED:LD set for $T$.
Based on the above analysis of $|V(T)|$ and $|S|$, we see that indeed $T \in \mathscr{T}_{min}^1$ for any choice of $\beta$.
\end{proof}

\begin{lemma}\label{lem:type-1-extremal-2}
Let $S$ be an optimal RED:LD set for a tree $T \in \mathscr{T}_{min}^1$ with order at least 7. Then every $v \notin S$ has $p = deg(v) \le 3$ and there exists $\beta =(\beta_1,\hdots,\beta_p)$ with $\beta \in \{(0,0),(1,2),(2,2,2)\}$ such that the connected components of $T-v$, $T_1, \hdots, T_p$, satisfy $T_i \in \mathscr{T}_{min}^i$, and $S_i = S \cap V(T_i)$ is optimal on $T_i$.
\end{lemma}
\begin{proof}
Since $T$ is in $\mathscr{T}_{min}^1$, let $|V(T)| = 3k +1$ with $|S| = 2k+2$, meaning there are $j = k\!-\!1$ non-detectors.
Let $v \notin S$ with $deg(v) = p$ and define $T_1,\hdots,T_p$, $j_1,\hdots,j_p$, and $S_1,\hdots,S_p$ similarly to Lemma~\ref{lem:type-2-extremal-2}.
We again find that $|V(T)| \ge 3(j - 1) + 2p +1$.
If $deg(v) \ge 4$, then  $|V(T)| \ge 3(j-1) + 2 \times 4 + 1 = 3k+3$, a contradiction.
If $deg(v)=3$, then $|V(T)|$ is exactly $3(j-1) + 2 \times 3 + 1 = 3k + 1$, so each of $T_i$ must have precisely $3j_i+2$ vertices (see Table~\ref{tab:redld-min-tree}), meaning $T_i \in \mathscr{T}_{min}^2$; this corresponds to $\beta = (2,2,2)$.
Otherwise $deg(v) = 2$, and because $T$ is assumed to be in $\mathscr{T}_{min}^1$, we know $|V(T)| = 3j + 4$.
Let $|V(T_i)| = 3j_i + \alpha_i$; Observation~\ref{obs:max-n-for-k-nondetector} gives us that $\alpha_i \ge 2$.
We know that $|V(T_1)| + |V(T_2)| + 1 = |V(T)|$, which means $(3j_1 + \alpha_1) + (3j_2 + \alpha_2) + 1 = 3(j_1 + j_2) + \alpha_1 + \alpha_2 + 1 = 3(j-1) + \alpha_1 + \alpha_2 + 1 = 3j + 4$.
Thus, we see that $\alpha_1 + \alpha_2 = 6$, so without loss of generality $(\alpha_1, \alpha_2) \in \{(2,4), (3,3)\}$; from Table~\ref{tab:redld-min-tree}, we see these correspond to the $\beta=(1,2)$ and $\beta=(0,0)$ cases.
In any case, we saw from the table that each $S_i$ must be optimal on $T_i$, completing the proof.
\end{proof}

From the previous six Lemmas, we have the following theorem.

\begin{theorem}\label{theo:t-min-family}
Starting with the four extremal trees on $n \le 4$ vertices given in Figure~\ref{fig:redld-extreme-low-tree}, the construction lemmas (\ref{lem:type-2-extremal-1}, \ref{lem:type-0-extremal-1}, and \ref{lem:type-1-extremal-1}) create the entire family of extremal trees, $\mathscr{T}_{min}$.
\end{theorem}
\begin{proof}
The aforementioned construction lemmas show that any tree they produce is extremal.
Lemmas \ref{lem:type-2-extremal-2}, \ref{lem:type-0-extremal-2}, and \ref{lem:type-1-extremal-2} show that any $T \in \mathscr{T}_{min}$ can be formed by using one of the construction lemmas.
\end{proof}

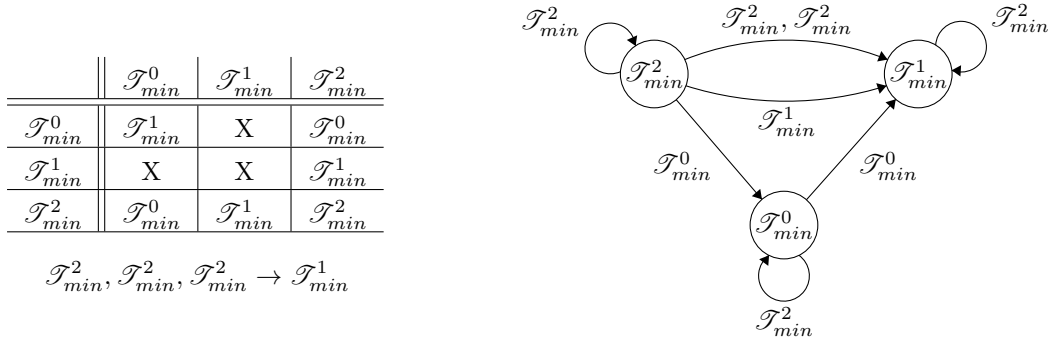
\begin{figure}[ht]
    \centering
    \begin{minipage}[t]{.25\textwidth}
        \vspace{20pt}
        \centering
        \begin{tabular}{c@{\rule{0pt}{1.2em}}@{\hspace{0.55em}}||c|c|c}
            & $\mathscr{T}_{min}^0$ & $\mathscr{T}_{min}^1$ & $\mathscr{T}_{min}^2$ \\ \hline \hline
            $\mathscr{T}_{min}^0$ & $\mathscr{T}_{min}^1$ & X & $\mathscr{T}_{min}^0$ \\ \hline
            $\mathscr{T}_{min}^1$ & X & X & $\mathscr{T}_{min}^1$ \\ \hline
            $\mathscr{T}_{min}^2$ & $\mathscr{T}_{min}^0$ & $\mathscr{T}_{min}^1$ & $\mathscr{T}_{min}^2$ \\ \hline
        \end{tabular}
        \begin{equation*}
            \;\;\;\;\;\mathscr{T}_{min}^2, \mathscr{T}_{min}^2, \mathscr{T}_{min}^2 \rightarrow \mathscr{T}_{min}^1
        \end{equation*}
    \end{minipage}
    \hspace{7em}
    \begin{minipage}[t]{.4\textwidth}
        \vspace{0pt}
        \centering
        \begin{tikzpicture}[scale=0.15]
            \tikzstyle{every node}+=[inner sep=0pt]
            \draw [black] (8.9,-5.8) circle (3);
            \draw (8.9,-5.8) node {$\mathscr{T}_{min}^2$};
            \draw [black] (32.3,-5.8) circle (3);
            \draw (32.3,-5.8) node {$\mathscr{T}_{min}^1$};
            \draw [black] (20.4,-19.2) circle (3);
            \draw (20.4,-19.2) node {$\mathscr{T}_{min}^0$};
            \draw [black] (29.489,-6.846) arc (-72.50103:-107.49897:29.564);
            \fill [black] (29.49,-6.85) -- (28.58,-6.61) -- (28.88,-7.56);
            \draw (20.6,-8.71) node [below] {$\mathscr{T}_{min}^1$};
            \draw [black] (33.847,-3.243) arc (176.55:-111.45:2.25);
            \draw (38.61,-0.86) node [right] {$\mathscr{T}_{min}^2$};
            \fill [black] (35.27,-5.47) -- (36.04,-6.02) -- (36.1,-5.02);
            \draw [black] (10.85,-8.08) -- (18.45,-16.92);
            \fill [black] (18.45,-16.92) -- (18.3,-15.99) -- (17.55,-16.64);
            \draw (14.1,-13.95) node [left] {$\mathscr{T}_{min}^0$};
            \draw [black] (22.39,-16.96) -- (30.31,-8.04);
            \fill [black] (30.31,-8.04) -- (29.4,-8.31) -- (30.15,-8.97);
            \draw (26.89,-13.95) node [right] {$\mathscr{T}_{min}^0$};
            \draw [black] (21.723,-21.88) arc (54:-234:2.25);
            \draw (20.4,-26.45) node [below] {$\mathscr{T}_{min}^2$};
            \fill [black] (19.08,-21.88) -- (18.2,-22.23) -- (19.01,-22.82);
            \draw [black] (5.919,-5.593) arc (293.76722:5.76722:2.25);
            \draw (2.42,-1.16) node [left] {$\mathscr{T}_{min}^2$};
            \fill [black] (7.25,-3.31) -- (7.39,-2.37) -- (6.47,-2.78);
            \draw [black] (11.611,-4.52) arc (111.73028:68.26972:24.278);
            \fill [black] (29.59,-4.52) -- (29.03,-3.76) -- (28.66,-4.69);
            \draw (20.6,-2.29) node [above] {$\mathscr{T}_{min}^2,\mathscr{T}_{min}^2$};
        \end{tikzpicture}
    \end{minipage}
    \caption{Construction of extremal trees---as per Theorem~\ref{theo:t-min-family}---shown in table and state machine form.}
    \label{fig:state-machine}
\end{figure}

The construction patterns for $\mathscr{T}_{min}$ proven in Theorem~\ref{theo:t-min-family} are summarized in Figure~\ref{fig:state-machine} in both table form and as a finite state machine, where the states represent the current tree type (mod size) and the edge labels are the types of tress being combined with it (through a new non-detector vertex).
For example, in row 2 column 3 of the table, we see that combining two trees in $\mathscr{T}_{min}^1$ and $\mathscr{T}_{min}^2$, respectively, with a non-detector yields a tree in $\mathscr{T}_{min}^1$; this is shown in the finite state machine as the $\mathscr{T}_{min}^2$ transition from $\mathscr{T}_{min}^1$ to $\mathscr{T}_{min}^1$, or as the $\mathscr{T}_{min}^1$ transition from $\mathscr{T}_{min}^2$ to $\mathscr{T}_{min}^1$.

\begin{theorem}\label{theo:all-solutions}
Let $S$ be an arbitrary optimal RED:LD set for a tree in $\mathscr{T}_{min}$. Then $S$ can be generated by Theorem~\ref{theo:t-min-family} using only RED:LD sets and trees found by the same Theorem for smaller trees.
\end{theorem}
\begin{proof}
The argument will proceed inductively; we will assume that the construction Theorem produces all possible optimal RED:LD sets for all extremal trees on up to $n-1$ vertices using only previously found RED:LD sets.
Clearly, if $n \le 4$ then we have all solutions, as these are the base case trees.
Otherwise, let $v \in V(T) - S$ be an arbitrary non-detector; then we can split the tree into $p=deg(v)$ subtrees, $T_1,\hdots,T_p$, of $T-v$; let $S_i = S \cap V(T_i)$.
The inductive hypothesis guarantees that the construction Theorem can produce all possible optimal RED:LD sets on each $T_i$, which means each $S_i$ can be formed by the construction Theorem using only previously found RED:LD sets.
The original tree, $T$, and RED:LD set, $S$, are simply formed by connecting all of the $T_i$ with a non-detector and unioning the $S_i$, so $T$ and $S$ can be constructed using only previously found trees and RED:LD sets, completing the proof.
\end{proof}

The findings of Theorem~\ref{theo:all-solutions} show that Theorem~\ref{theo:t-min-family} not only produces the family of all extremal trees, but also produces every optimal RED:LD set on each tree (without ever needing an arbitrary optimal RED:LD set, as was used in proving the original three construction Lemmas).

Theorem~\ref{theo:t-min-family} produces extremal trees and RED:LD sets by connecting smaller extremal trees with a new non-detector vertex.
From Theorem~\ref{theo:all-solutions}, we see that this process can be done using only previously found RED:LD sets on the smaller trees.
Thus, detectors in the graph must come from the original four base case trees, giving us the following corollary:

\begin{corollary}
If $S$ is an optimal RED:LD set for $T \in \mathscr{T}_{min}$, then each connected component of the graph generated by $S$ is one of the four trees given in Figure~\ref{fig:redld-extreme-low-tree}.
\end{corollary}

\begin{corollary}
If $T \in \mathscr{T}_{min}^2$, then any minimum RED:LD set efficiently 2-dominates all vertices.
\end{corollary}


\begin{wrapfigure}[7]{r}{0.16\textwidth}
    \centering
    \includegraphics[width=0.14\textwidth]{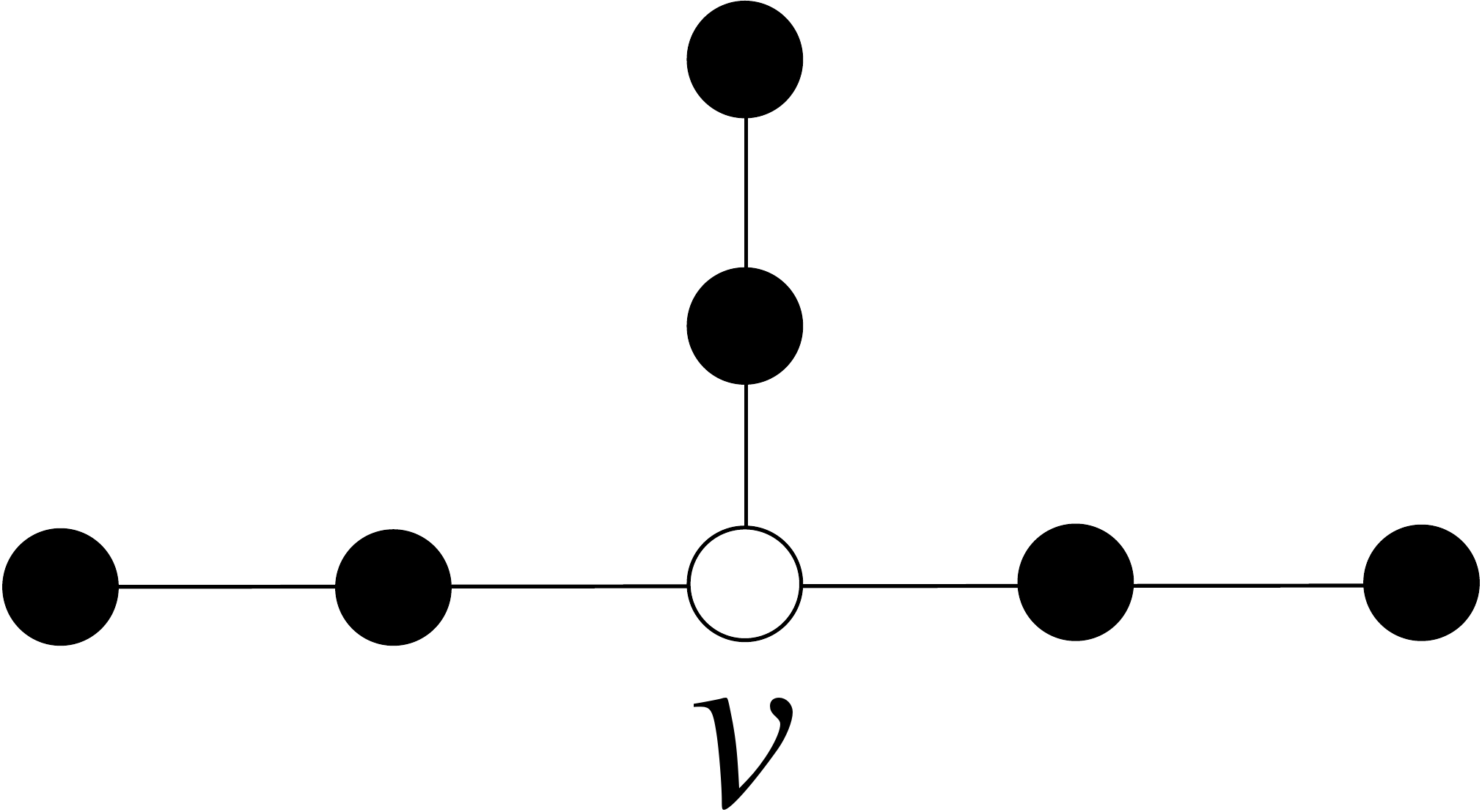}
    \caption{$T_7$.}
    \label{fig:T7}
\end{wrapfigure}

Figure~\ref{fig:redld-min-trees-all} shows all trees in $\mathscr{T}_{min}$ on 5--10 vertices (trees on 2--4 vertices are given by the base case trees from Figure~\ref{fig:redld-extreme-low-tree}).
Figure~\ref{fig:T7} is a particular labeled subgraph that is only present for extremal trees constructed using the rule for combining three trees in $\mathscr{T}_{min}^2$ with a single non-detector.

\begin{corollary}
For any optimal RED:LD set on a tree $T \in \mathscr{T}_{min}$, if $v$ is a degree three non-detector then $T \in \mathscr{T}_{min}^1$, $v$ is unique, and the three subtrees of $T-v$ can be efficiently 2-dominated.
\end{corollary}

We will now give algorithms for checking if an arbitrary tree, $T$, is in $\mathscr{T}_{min}$ or not.
This is given by Algorithms \ref{alg:redld-tree-type-2}, \ref{alg:redld-tree-type-0}, and \ref{alg:redld-tree-type-1}, one for each of $\mathscr{T}_{min}^2$, $\mathscr{T}_{min}^0$, and $\mathscr{T}_{min}^1$, respectively.
Each of these algorithms returns a (non-empty) vertex set which is an optimal RED:LD set on the tree if $T \in \mathscr{T}_{min}$, or an empty set if $T \notin \mathscr{T}_{min}$.
Algorithms \ref{alg:redld-tree-type-2}--\ref{alg:redld-tree-type-1} all begin with similar preprocessing logic, which is performed by Algorithm~\ref{alg:redld-tree-helper}; this finds and removes all exterior $P_2$ trees that were connected by a $\mathscr{T}_{min}^2$ transition (see the finite state machine from Figure~\ref{fig:state-machine}), along with the non-detector connecting them.

\begin{algorithm}[p]
\caption{Get the set of all recursive exterior $P_2$ and non-detector components in a tree}
\label{alg:redld-tree-helper}
\begin{algorithmic}[1]
\Function{ExteriorP2PlusNondetector}{$T, q$}
\State{\textbf{if} $|V(T)| \ge q$ and $\exists u,v,w$, where $N(v)=\{u,w\}$, $deg(u)=1$, and $deg_0(w) = 2$}
    \State{\hspace{2em}\textbf{let} $T' = T-\{u,v,w\}$}
    \State{\hspace{2em}\textbf{let} $S_1,S_2 = \textproc{ExteriorP2PlusNondetector}(T', q)$}
    \State{\hspace{2em}\textbf{return} $(S_1 \cup \{u,v\}, S_2 \cup \{w\})$}
\State {\textbf{else return} $(\varnothing,\varnothing)$}
\EndFunction
\end{algorithmic}
\end{algorithm}

\begin{algorithm}[p]
\caption{Classify a tree on $3k+2$ vertices as extremal or not extremal}
\label{alg:redld-tree-type-2}
\begin{algorithmic}[1]
\Function{ExtremalRedLd2}{$T$}
\State{\textbf{let} $S_1,S_2 = \textproc{ExteriorP2PlusNondetector}(T, 0)$}
\State{$T \gets T-(S_1 \cup S_2)$}
\State{\textbf{if} $|V(T)| \le 1$ \textbf{then return} $\varnothing$}
\State{\textbf{if} $|V(T)| = 2$ \textbf{then return $V(T) \cup S_1$}}
\EndFunction
\end{algorithmic}
\end{algorithm}

\begin{algorithm}[p]
\caption{Classify a tree on $3k$ vertices as extremal or not extremal}
\label{alg:redld-tree-type-0}
\begin{algorithmic}[1]
\Function{ExtremalRedLd0}{$T$}
\State{\textbf{let} $S_1,S_2 = \textproc{ExteriorP2PlusNondetector}(T, 0)$}
\State{$T \gets T-(S_1 \cup S_2)$}
\State{\textbf{if} $|V(T)| \le 2$ \textbf{then return} $\varnothing$}
\State{\textbf{if} $|V(T)| = 3$ \textbf{then return $V(T) \cup S_1$}}
\EndFunction
\end{algorithmic}
\end{algorithm}

\begin{algorithm}[p]
\caption{Classify a tree on $3k+1$ vertices as extremal or not extremal}
\label{alg:redld-tree-type-1}
\begin{algorithmic}[1]
\Function{ExtremalRedLd1}{$T$}
\State{\textbf{let} $S_1,S_2 = \textproc{ExteriorP2PlusNondetector}(T, 5)$}
\State{$T \gets T-(S_1 \cup S_2)$}
\State{\textbf{if} $|V(T)| \le 3$ \textbf{then return} $\varnothing$}
\State{\textbf{if} $|V(T)| = 4$ \textbf{then return $V(T) \cup S_1$}}
\State{\textbf{if} $T = T_7$ \textbf{then return $(V(T) - \{v\}) \cup S_1$}}
\State{\textbf{if} $\exists B_1,B_2 \subseteq V(T)$ be two branches each on 3 vertices with parent $w$ having $deg_0(w)=2$}
\State{\hspace{2em}\textbf{if} $w_1=w_2$ \textbf{return} $\{u_1,v_1,x_1,u_2,v_2,x_2\} \cup S_1$}
\State{\hspace{2em}\textbf{let} $T' = T-(B_1 \cup B_2 \cup \{w_1, w_2\})$}
\State{\hspace{2em}\textbf{let} $S' = \textproc{ExtremalRedLd2}(T')$}
\State{\hspace{2em}\textbf{if} $S' = \varnothing$ \textbf{return} $\varnothing$}
\State{\hspace{2em}\textbf{else return} $\{u_1,v_1,x_1,u_2,v_2,x_2\} \cup S' \cup S_1$}
\State{\textbf{else return} $\varnothing$}
\EndFunction
\end{algorithmic}
\end{algorithm}

\subsubsection{Infinite k-aray trees}

\begin{figure}[ht]
    \centering
    \includegraphics[width=0.9\textwidth]{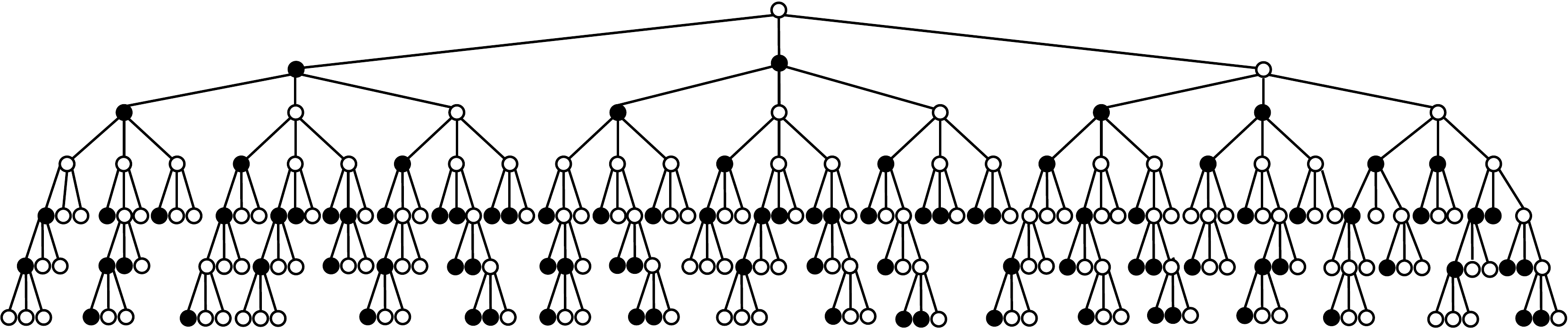}
    \caption{A partial view of the infinite, complete $k$-ary tree with $k = 3$.}
    \label{fig:infinite-k-ary}
\end{figure}

\begin{theorem}\label{theo:inf-trees}
Let $T$ be an infinite, complete $k$-ary tree for $k \ge 2$; then $\textrm{RED:LD}(T) = \frac{2}{k+2}$.
\end{theorem}
\begin{proof}
As any RED:LD set must be 2-dominating, for any detector vertex $x$, $sh(x) \le (k+2) \times \frac{1}{2}$, so we find that $\textrm{RED:LD}(T) \ge \frac{2}{k+2}$.

Let $r$ denote the root vertex, which is the unique vertex with degree $k$, for an infinite complete $k$-ary tree.
We construct a set $S \subseteq V(T)$ by first deciding $r \notin S$ and $|N(r) \cap S| = 2$; we say that all vertices in $N[r]$ have been ``visited".
From here we apply a recursive pattern: where $x$ is a vertex (already visited) and $C$ is its set of children (not yet visited), let $|N[x] \cap S| = 2$, with any needed detectors coming from $C$; we then mark all vertices in $C$ as visited.
This construction is demonstrated in Figure~\ref{fig:infinite-k-ary} for $k=3$, and clearly, every vertex is at least 2-dominated.
Suppose $u,v \notin S$; then by this construction $\exists x \in (N(u)-N(v)) \cap S$ and $\exists y \in (N(v)-N(u)) \cap S$, so $u$ and $v$ are 2-distinguished.
Otherwise, suppose $u \in S$ and $v \notin S$; then it must be the case that $\exists x \in (N(v)-N(u)) \cap S$, so $u$ and $v$ are 1-distinguished.
Thus, by Theorem~\ref{theo:red-ld-char}, $S$ is a RED:LD set.
Since every detector in $S$ has share $(k+2) \times \frac{1}{2}$ in this construction, the density of $S$ in $V(T)$ clearly is $\frac{2}{k+2}$; so we have $\textrm{RED:LD}(T) \le \frac{2}{k+2}$.
\end{proof}

\begin{figure}[ht]
    \centering
    \begin{tabular}{cc}
        \includegraphics[width=0.7\textwidth]{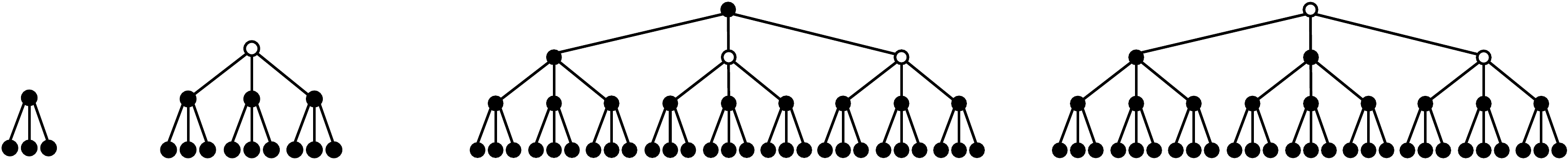} \\ \\
        \includegraphics[width=0.7\textwidth]{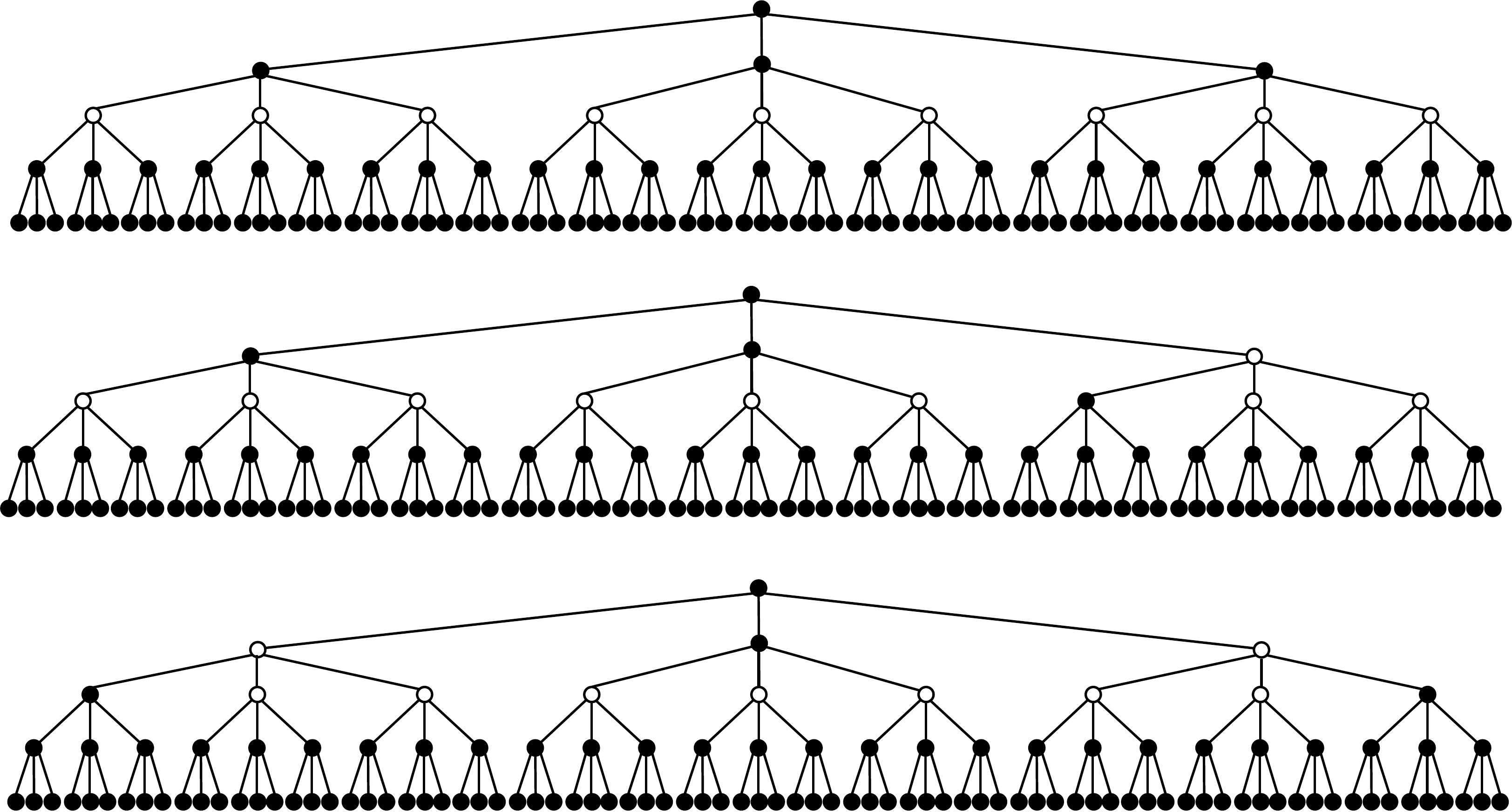}
    \end{tabular}
    \caption{Exhaustive optimal RED:LD sets on complete 3-ary trees with depths 1--4.}
    \label{fig:finite-k-ary}
\end{figure}

\FloatBarrier
\begin{theorem}\label{theo:kary-finite}
Let $T$ be a complete $k$-ary tree for $k \ge 2$ with depth $d \ge 1$.
Then, where $m = d\!\!\mod 3$ and $t = (d\!-\!1)\!\!\mod 3$,
\begin{equation*}
    \textrm{RED:LD}(T) = 
    \begin{cases}
    k+1 & d = 1 \\
    k^2 + k & d = 2 \\
    k^3 + k^2 + 2 & d = 3 \\
    (1\!-\!m)(2\!-\!m) + k^t(k\!+\!1)\frac{k^{3\ceil{d/3}}-1}{k^3-1} & d \ge 4
    \end{cases}
\end{equation*}
\end{theorem}
\begin{proof}
Let $S$ be an optimal RED:LD set on $T$.
If $d = 1$ then the results are clear from the requirement of 2-domination alone.
If $d \ge 2$ then where $L = \{v \in V(T) : deg(v) = 1\}$, we require $\cup_{v \in L}{N[v]} \subseteq S$ to 2-dominate each vertex in $L$; this includes every vertex at depths $d$ and $d-1$, so $|S| \ge k^d + k^{d-1}$.
If $d = 2$, it is clear that this set satisfies Proposition~\ref{prop:tree-2-dom}, so $\textrm{RED:LD}(T) = k^2 + k$.
If $d = 3$, then we will still need at least two additional detectors to 2-dominate the root vertex; this is sufficient to 2-dominate every vertex, so $\textrm{RED:LD}(T) = k^3 + k^2 + 2$.
Figure~\ref{fig:finite-k-ary} shows all non-isomorphic solutions for minimum RED:LD sets on a 3-ary tree with depth $d \le 4$.

Otherwise, we assume $d \ge 4$.
Consider a vertex $x$ at depth $d-4$, with $T_x$ being the sub-tree with $x$ as its root; we will show that any optimal $S$ requires $x \in S$.
As stated previously, to 2-dominate the leaf vertices we have all vertices at depths $d$ and $d-1$ being detectors.
Let $C$ denote the set of $k$ children of $x$; clearly for any $v \in C$, we require $|N[v] \cap S| \ge 2$.
If $x \in S$ then we need at least $k+1$ additional detectors, including $x$ itself; the last three graphs of Figure~\ref{fig:finite-k-ary} show all non-isomorphic ways to use only $k+1$ additional detectors in $T_x$.
By having $\{x\} \cup C \subseteq S$, we see that all vertices in $T_x$ are 2-dominated regardless of vertices outside of $T_x$, if any; thus $k+1$ additional detectors are sufficient to be a RED:LD set.
Otherwise, $x \notin S$, so we require at least $2k$ additional detectors.
Because $k \ge 2$ by hypothesis, $k+1 < 2k$, and having $x \notin S$ potentially means we require even more detectors outside of $T_x$.
Thus, if $x \notin S$, then $S$ is not optimal, a contradiction.

Therefore, we assume $x \in S$ and use exactly $k+1$ additional detectors, including $x$ itself.
As previously stated, taking taking $\{x\} \cup C \subseteq S$ is sufficient to 2-dominate all vertices in $T_x$ regardless of vertices outside of $T_x$; applying this across all $x$ vertices, we see that all vertices at depths $d-4$ and $d-3$ can be detectors while still having $S$ be optimal.
We see that rows $d-4$ and $d-3$ are similar to rows $d$ and $d-1$, so we can repeat the logic of this proof starting at the beginning for the tree of depth $d-3$, giving us the recurrence relation $\textrm{RED:LD}(T) = \textrm{RED:LD}(T') + k^d + k^{d-1}$ for $d \ge 4$ where $T'$ is the truncated tree with root $r$ and depth $d - 3$.
Let $f(d) = \textrm{RED:LD}(T)$ with $d$ denoting depth; then we can expand the recurrence relation as

\begin{align*}
    f(d) &= f(d-3) + k^d + k^{d-1} \\
    &= f(d-6) + k^{d-3} + k^{d-4} + k^d + k^{d-1} \\
    &= f(d-3p) + k^{d-3p+3} + k^{d-3p+2} + \cdots + k^d + k^{d-1}
\end{align*}

We will now consider the cases of $d-3p$ being 1, 2, or 3 for some $p \in \mathbb{N}$, for which we can apply the previous expansion of $f(d)$ using the base cases $f(1)$, $f(2)$, and $f(3)$ found above.

\begin{align*}
    d-3p = 1\!:\;\;\;\; f(d) &= f(d-3p) + k^{d-3p+3} + k^{d-3p+2} + \cdots + k^d + k^{d-1} \\
    &= (k+1) + k^4 + k^3 + k^7 + k^6 + \cdots + k^{3p+1} + k^{3p} \\
    &= (k+1) + (k^4 + k^7 + k^{10} + \cdots + k^{3p+1}) + (k^3 + k^6 + k^9 + \cdots + k^{3p}) \\
    &= (k+1) + (k + 1)(k^3 + k^6 + \cdots + k^{3p} + 1 - 1) \\
    &= (k+1) + (k + 1)\left(\frac{k^{3p+3}-1}{k^3-1} - 1\right) \\
    &= (k + 1)\frac{k^{d+2}-1}{k^3-1} \\
    d-3p = 2\!:\;\;\;\; f(d) &= (k^2 + k) + k^5 + k^4 + k^8 + k^7 + \cdots + k^{3p+2} + k^{3p+1} \\
    &= k(k+1) + k(k + 1)\left(\frac{k^{3p+3}-1}{k^3-1} - 1\right) \\
    &= k(k + 1)\frac{k^{d+1}-1}{k^3-1} \\
    d-3p = 3\!:\;\;\;\; f(d) &= (k^3 + k^2 + 2) + k^6 + k^5 + k^9 + k^8 + \cdots + k^{3p+3} + k^{3p+2} \\
    &= 2 + k^2(k + 1) + k^2(k + 1)\left(\frac{k^{3p+3}-1}{k^3-1} - 1\right) \\
    &= 2 + k^2(k + 1)\frac{k^{d}-1}{k^3-1}
\end{align*}

The three closed forms found above can be combined into one expression using modulo-3 arithmetic, giving the form of the $d \ge 4$ case of the theorem statement.
Table~\ref{tab:kary-finite} gives example values for $\textrm{RED:LD}(T)$ for various choices of $k$ and $d$.
\end{proof}

\begin{table}[ht]
    \centering
    \begin{tabular}{r|rrrrrr}
        $d$ & $k=2$ & $k=3$ & $k=4$ & $k=5$ & $k=6$ & $k=7$ \\ \hline
        1 & 3 (0.75) & 4 (0.89) & 5 (0.94) & 6 (0.96) & 7 (0.97) & 8 (0.98) \\
        2 & 6 (0.75) & 12 (0.89) & 20 (0.94) & 30 (0.96) & 42 (0.97) & 56 (0.98) \\
        3 & 14 (0.88) & 38 (0.94) & 82 (0.96) & 152 (0.97) & 254 (0.98) & 394 (0.98) \\
        4 & 27 (0.84) & 112 (0.92) & 325 (0.95) & 756 (0.97) & 1519 (0.98) & 2752 (0.98) \\
        5 & 54 (0.84) & 336 (0.92) & 1300 (0.95) & 3780 (0.97) & 9114 (0.98) & 19264 (0.98) \\
        6 & 110 (0.86) & 1010 (0.92) & 5202 (0.95) & 18902 (0.97) & 54686 (0.98) & 134850 (0.98) \\
        7 & 219 (0.86) & 3028 (0.92) & 20805 (0.95) & 94506 (0.97) & 328111 (0.98) & 943944 (0.98) \\
        8 & 438 (0.86) & 9084 (0.92) & 83220 (0.95) & 472530 (0.97) & 1968666 (0.98) & 6607608 (0.98) \\
        9 & 878 (0.86) & 27254 (0.92) & 332882 (0.95) & 2362652 (0.97) & 11811998 (0.98) & 46253258 (0.98) \\
        10 & 1755 (0.86) & 81760 (0.92) & 1331525 (0.95) & 11813256 (0.97) & 70871983 (0.98) & 323772800 (0.98)
    \end{tabular}
    \caption{Sample values from Theorem~\ref{theo:kary-finite}. Entries are the exact value of $\textrm{RED:LD}(T_d)$ and the density.}
    \label{tab:kary-finite}
\end{table}

\FloatBarrier
\subsection{Infinite grids}\label{sec:inf-grids}

\begin{figure}[ht]
    \centering
    \begin{tabular}{c@{\hskip 5em}c}
        \includegraphics[width=0.27675\textwidth]{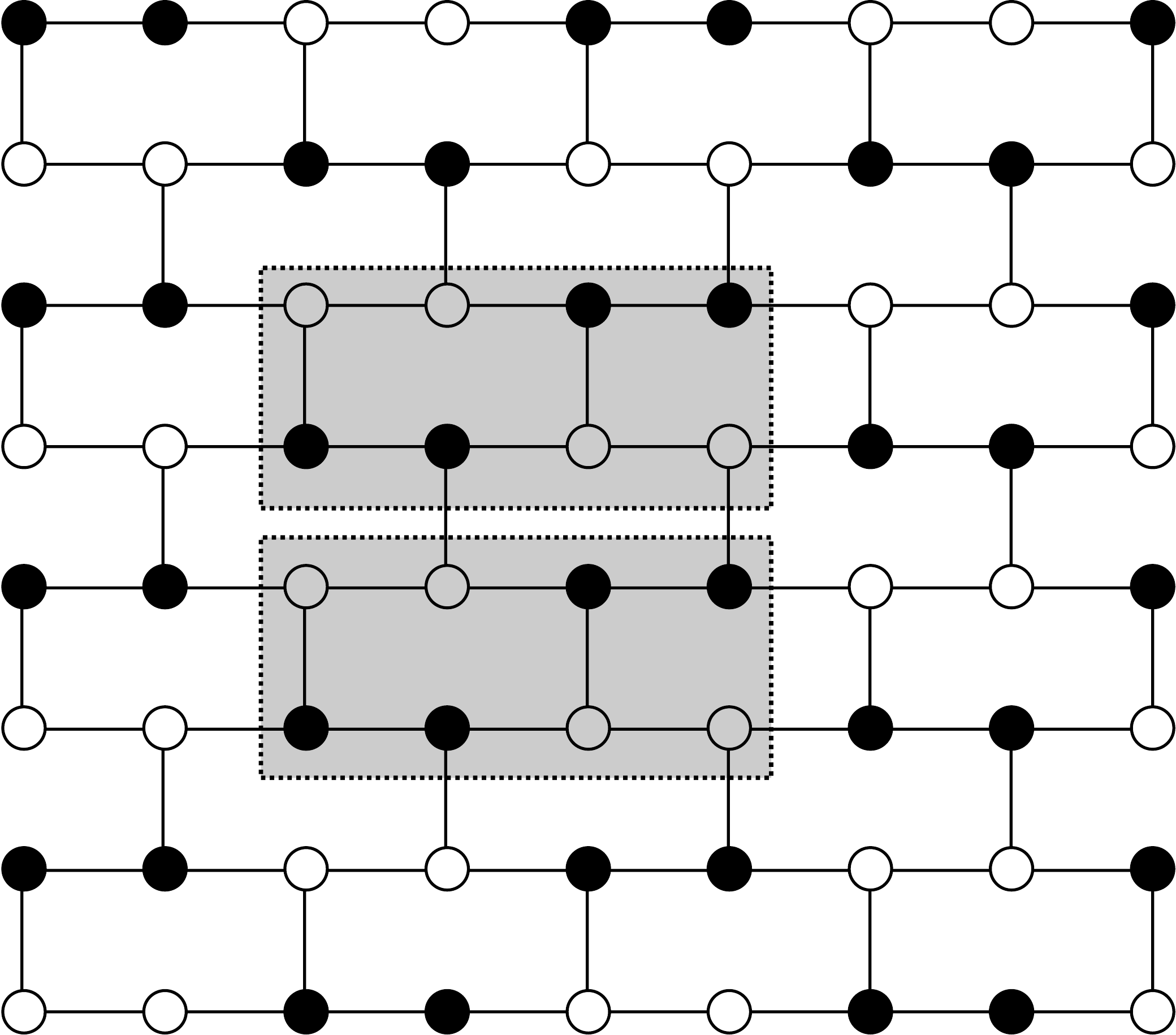} & \includegraphics[width=0.35\textwidth]{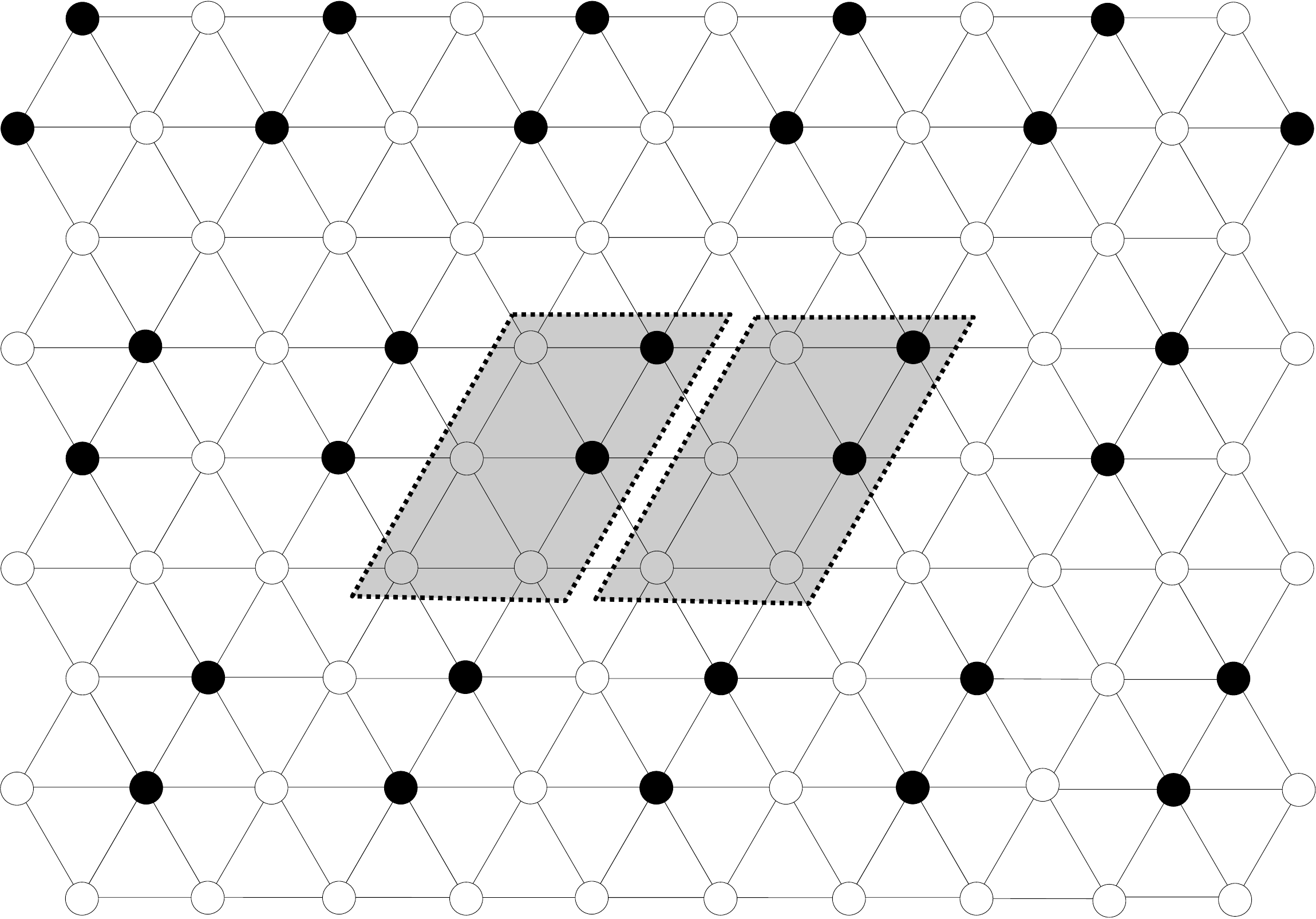} \\ (a) & (b) \\ \\
        \includegraphics[width=0.35\textwidth]{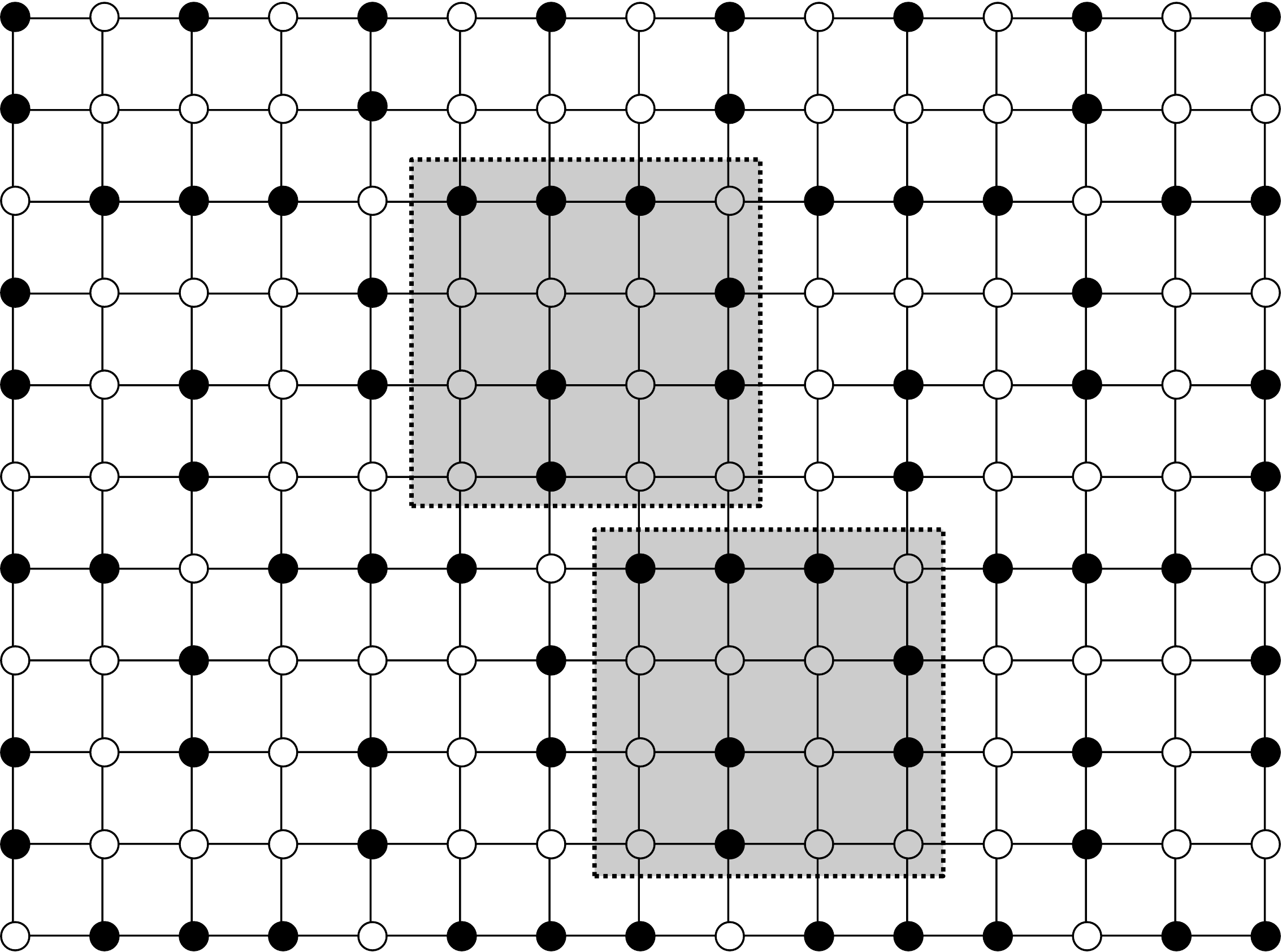} &  \includegraphics[width=0.35\textwidth]{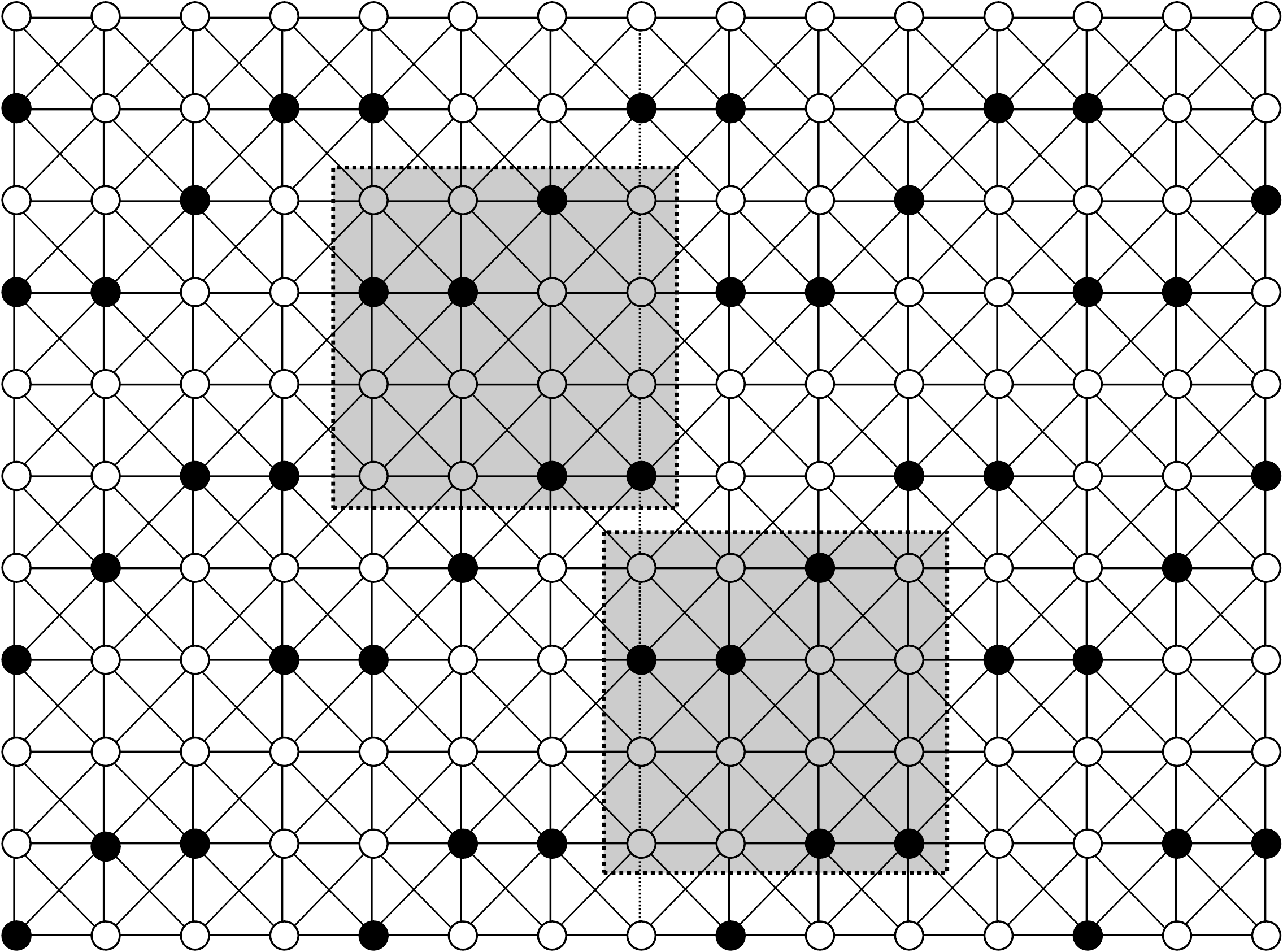} \\ (c) & (d)
    \end{tabular}
    \caption{Upper bounds for the RED:LD on (a) $HEX$, (b) $TRI$, (c) $SQ$, and (d) $K$.}
    \label{fig:inf-grids-upper-bounds}
\end{figure}

\begin{theorem}
For the infinite hexagonal grid, $HEX$, $\textrm{RED:LD}(HEX) = \frac{1}{2}$.
\end{theorem}
\cbeginproof
From Figure~\ref{fig:inf-grids-upper-bounds}~(a), we have that $\textrm{RED:LD}(HEX) \le \frac{1}{2}$; we need only show that $\textrm{RED:LD}(HEX) \ge \frac{1}{2}$.
Any RED:LD set must be 2-dominating, and the HEX grid is 3-regular.
Thus, the share of any detector vertex is at most $4 \times \frac{1}{2} = 2$, giving a lower bound of $\frac{1}{2}$.
\cendproof

\begin{theorem}
For the infinite triangular grid, $TRI$, $\textrm{RED:LD}(TRI) = \frac{1}{3}$.
\end{theorem}
\begin{wrapfigure}{r}{0.25\textwidth}
    \centering
    \includegraphics[width=0.25\textwidth]{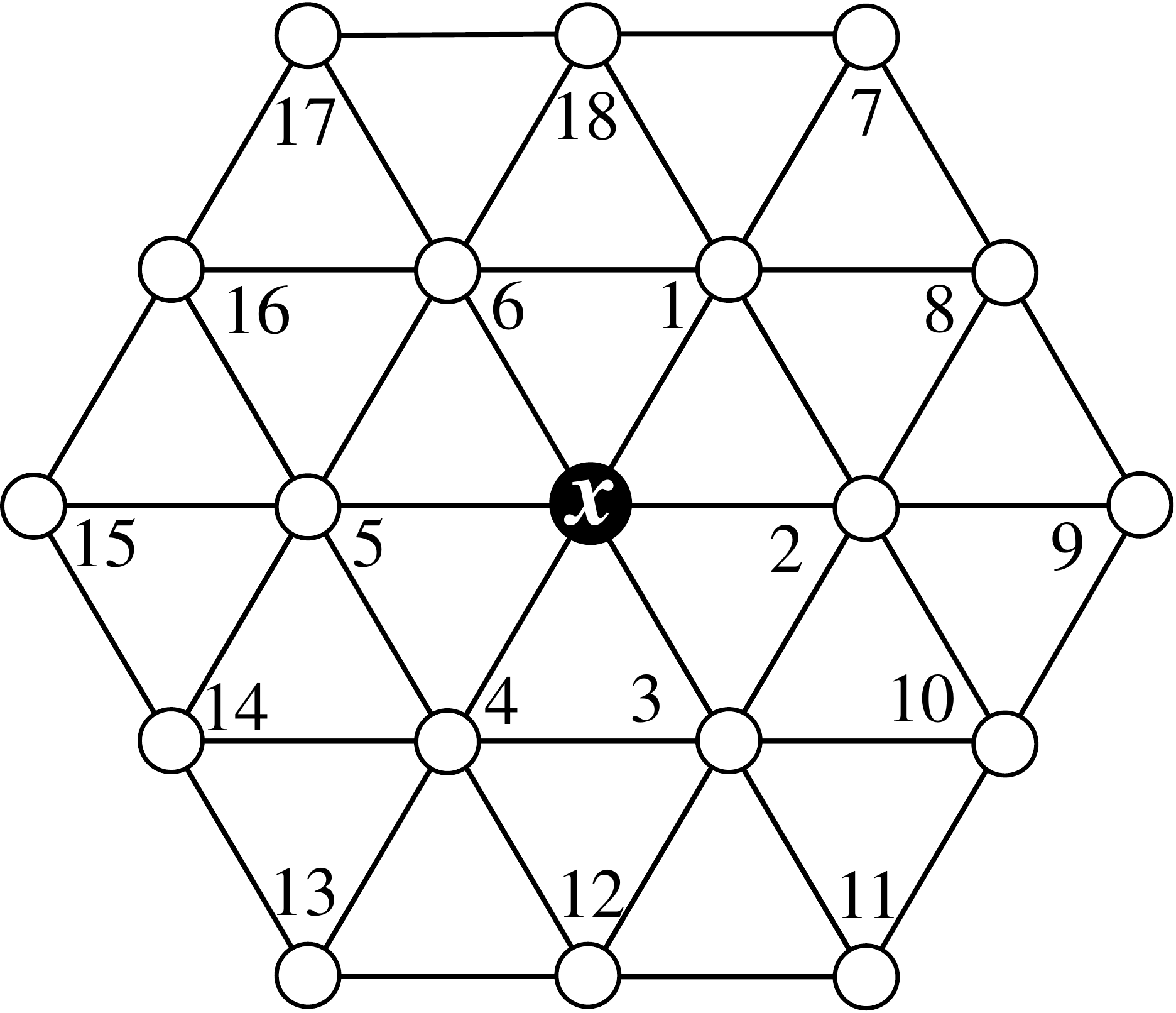}
    \caption{Vertex labeling.}
    \label{fig:tri-proof-labeling}
\end{wrapfigure}
\cbeginproof
From Figure~\ref{fig:inf-grids-upper-bounds}~(b), we see that $\textrm{RED:LD}(TRI) \le \frac{1}{3}$.
Thus, we need only show that $\textrm{RED:LD}(TRI) \ge \frac{1}{3}$; to do this, we will demonstrate that the average share of all detectors in a RED:LD set $S$ on $TRI$ is at most $3$.
Let $x \in S$; we continue by considering the possible values of $dom(x)$.
Refer to Figure~\ref{fig:tri-proof-labeling}; note that $v_k$ will refer to the vertex labeled $k$.
First, suppose $dom(x) \ge 4$; then we have $|N(x) \cap D_{3+}| \ge 2$, so $sh(x) \le \sigma_{4332222} < 3$ and we are done.

Next, suppose $dom(x) = 3$; then there are three non-isomorphic cases.
\textbf{Case~1:} $\{v_1,v_2\} \subseteq S$.
We see that $dom(v_1) \ge 3$ and $dom(v_2) \ge 3$, so $sh(x) \le \sigma_{3332222} = 3$, and we are done.
\textbf{Case~2:} $\{v_1,v_3\} \subseteq S$.
We see that $v_1$ and $v_6$ are not distinguished; by property~2 of Theorem~\ref{theo:red-ld-char}, we must have $dom(v_1) \ge 3$ or $dom(v_6) \ge 3$; by symmetry, $dom(v_3) \ge 3$ or $dom(v_4) \ge 3$ as well.
Thus, $sh(x) \le \sigma_{3332222} = 3$, and we are done.
\textbf{Case~3:} $\{v_1,v_4\} \subseteq S$.
We see that $v_2$ and $v_6$ are not distinguished; by property~3 of Theorem~\ref{theo:red-ld-char}, we must have $sh[v_2v_6] \le \max\{\sigma_{33},\sigma_{24}\} = \sigma_{24}$; by symmetry, $sh[v_3v_5] \le \sigma_{24}$ as well.
Therefore, $sh(x) \le \sigma_{322} + \sigma_{24} + \sigma_{24} = \frac{17}{16} < 3$, and we are done.

Lastly, suppose $dom(x) = 2$; then we can assume $v_1 \in S$.
We see that to distinguish $x$ and $v_2$ we require $dom(v_2) \ge 3$; by symmetry, $dom(v_6) \ge 3$ as well.
If $\{v_1,v_3,v_4,v_5\} \cap D_{3+} \neq \varnothing$ then $sh(x) \le \sigma_{3332222} = 3$, and we are done; thus, we assume $\{v_1,v_3,v_4,v_5\} \subseteq D_{2}$.
Vertex $v_1$ is already 2-dominated, so $\{v_7,v_8,v_{18}\} \cap S = \varnothing$.
If $v_{12} \in S$ then $v_3$ and $v_4$ cannot be distinguished, a contradiction; thus, we assume $v_{12} \notin S$, and by symmetry $v_{14} \notin S$.
We require $v_{13} \in S$ to 2-dominate $v_4$.
If $\{v_{10},v_{16}\} \subseteq S$ then $dom(v_2) \ge 4$ to distinguish $v_2$ and $v_3$, and $dom(v_6) \ge 4$ to distinguish $v_5$ and $v_6$; we see that $sh(x) \le \sigma_{2442222} = 3$ and we are done.
Thus, without loss of generality, we can assume $v_{16} \notin S$, imposing $v_{15} \in S$ to 2-dominate $v_5$ and $v_{17} \in S$ to 3-dominate $v_6$.
We see that $dom(v_{18}) \ge 3$ to distinguish $v_6$ and $v_{18}$.
Suppose $v_{10} \notin S$; then $\{v_9,v_{11}\} \subseteq S$ to dominate $v_2$ and $v_3$.
We see that to distinguish $v_2$ and $v_8$ we require $dom(v_8) \ge 3$.
Thus, $sh(x) = \sigma_{2332222} = \frac{19}{6}$ and $sh(v_1) \le \sigma_{2333322} = \frac{17}{6}$.
As $v_1$ has only one neighboring detector, we may average their share values, which yields $\frac{1}{2}\left[ \frac{17}{6} + \frac{19}{6} \right] = 3$, and we are done.
Otherwise, we assume $v_{10} \in S$; then $v_{11} \notin S$ and $v_9 \in S$ to distinguish $v_2$ and $v_3$.
Thus, $sh(x) = \sigma_{2432222} = \frac{37}{12}$ and $sh(v_1) \le \sigma_{2433222} = \frac{35}{12}$.
Again, we find that it is safe to perform averaging, which yields $\frac{1}{2}\left[ \frac{37}{12} + \frac{35}{12} \right] = 3$, completing the proof.
\cendproof

\begin{theorem}
For the infinite square grid, $SQ$, $\frac{2}{5} \le \textrm{RED:LD}(SQ) \le \frac{7}{16}$.
\end{theorem}
\begin{proof}
From Figure~\ref{fig:inf-grids-upper-bounds}~(c) we see a RED:LD set on $SQ$ with density $\frac{7}{16}$, so $\textrm{RED:LD}(SQ) \le \frac{7}{16}$.
For any RED:LD set $S$, each vertex must be at least 2-dominated, and $SQ$ is 4-regular; thus, $\forall x \in S$, $sh(x) \le 5 \times \frac{1}{2}$.
Therefore, $\frac{2}{5} \le \textrm{RED:LD}(SQ)$, completing the proof.
\end{proof}

\begin{theorem}[\cite{ourking}]
For the infinite king grid, $K$, $\frac{3}{11} \le \textrm{RED:LD}(K) \le \frac{5}{16}$.
\end{theorem}

\bibliographystyle{acm}
\bibliography{refs}

\end{document}